\documentclass[a4paper,12pt]{article}
\usepackage{arxiv}

\usepackage[utf8]{inputenc} 
\usepackage[T1]{fontenc}    
\usepackage{hyperref}       
\usepackage{url}            
\usepackage{amsfonts}       
\usepackage{amsmath}        
\usepackage{graphicx}       
\usepackage{caption}
\usepackage{subcaption}
\usepackage{xcolor}

\definecolor{linkclr}{cmyk}{0.973,0.957,0,0.04}

\usepackage{environ}
\newcommand{\acksection}{\section*{Acknowledgments}}
\NewEnviron{ack}{
  \acksection
  \BODY
}

\title{Can Physics-Informed Neural Networks beat the Finite Element Method?} 
\author{Tamara G. Grossmann\thanks{Department of Applied Mathematics and Theoretical Physics, University of Cambridge, Cambridge, UK \hspace*{1.8em} \texttt{tg410@cam.ac.uk}} , Urszula Julia Komorowska\thanks{Department of Computer Science and Technology, University of Cambridge, Cambridge, UK} , Jonas Latz\thanks{Maxwell Institute for Mathematical Sciences and School of Mathematical and Computer Sciences, Heriot-Watt University, Edinburgh, UK} \hspace{1pt} and Carola-Bibiane Schönlieb\footnotemark[1]}

\usepackage{url}
\usepackage{hyperref}
\hypersetup{%
        colorlinks,
        unicode=false,
        breaklinks=true,
        plainpages=false,%
        citecolor=linkclr,
        linkcolor=linkclr,
        urlcolor=linkclr}

\usepackage{listings}
\begin{document}
\maketitle

\vspace{20pt}
\begin{abstract}
Partial differential equations play a fundamental role in the mathematical modelling of many processes and systems in physical, biological and other sciences. To simulate such processes and systems, the solutions of PDEs often need to be approximated numerically. The finite element method, for instance, is a usual standard methodology to do so. The recent success of deep neural networks at various approximation tasks has motivated their use in the numerical solution of PDEs. These so-called physics-informed neural networks and their variants have shown to be able to successfully approximate a large range of partial differential equations. So far, physics-informed neural networks and the finite element method have mainly been studied in isolation of each other. In this work, we compare the methodologies in a systematic computational study. Indeed, we employ both methods to numerically solve various linear and nonlinear partial differential equations: Poisson in 1D, 2D, and 3D, Allen--Cahn in 1D, semilinear Schr\"odinger in 1D and 2D.  We then compare computational costs and approximation accuracies. In terms of solution time and accuracy, physics-informed neural networks have not been able to outperform the finite element method in our study. In some experiments, they were faster at evaluating the solved PDE.
\end{abstract}

\section{Introduction}
\emph{Partial differential equations} (PDEs) are a corner stone of mathematical modelling and possibly applied mathematics itself. They are used to model a multitude of physical \cite{Lin1988}, biological \cite{taubes_2008}, socioeconomic \cite{Burger2014-qy}, and financial \cite{Heston} systems and processes. Beyond classical modelling, PDEs can describe the evolution of certain stochastic processes \cite{Risken1996}, be used in image reconstruction \cite{RUDIN1992}, as well as for filtering \cite{Kushner} and optimal control \cite{Bellmann} of dynamical systems.

The underlying idea is always the same: we aim to represent some process through a function that describes the behaviour of the process in space and time. The PDE is then a collection of laws that this function is supposed to satisfy. An obvious questions is, whether these laws are actually sufficient to uniquely specify this function \cite{Evans}.  Once uniqueness or even well-posedness \cite{Hadamard} has been established, the question is how to find the function satisfying these laws that will ultimately be the mathematical model.  

In many practical situations, it is impossible to find closed-form solutions for arising partial differential equations; they need to be solved numerically. 
Throughout the last decades, several numerical methods have been proposed and analysed to solve PDEs, especially the \emph{finite element method}  (FEM) \cite{Courant,Hrennikoff} that may be the standard methodology for a huge class of PDEs. Other techniques are, e.g., the finite difference \cite{iserles_2008}, finite volume \cite{Eymard}, and the spectral element \cite{PATERA1984} method. Due to the simplicity of the underlying approximation, many of these approaches are well-understood from a theoretical perspective: there are existing error estimators, as well as convergence and stability guarantees. Moreover, the given discretised problems are often in a form that can easily be solved numerically: relying on large, but sparse linear systems or on Newton solves with good initial values and convergence guarantees. Whilst implementing an FEM solver from scratch can be fairly tedious, several multipurpose computational libraries have appeared throughout the years, such as FEniCS \cite{Alnaes2015} or DUNE \cite{Sander2020}.

A clear disadvantage of this classical methodology is that it usually relies on a spatial discretisation through, e.g., a spatial grid or a large polynomial basis, and thus, let it suffer from the \emph{curse of dimensionality}: in three space dimensions, it can already be difficult to employ, e.g., the finite element method. In filtering and optimal control, we are easily interested in PDEs occurring in hundreds or thousands of dimensions -- here, classical methodology can rarely be employed. In addition, certain non-linear and non-smooth PDEs are difficult to discretise with, e.g., finite elements due to behaviour that needs to be resolved on a very fine grid, general non-smooth behaviour, or singularities. Since FEM is a mesh-based solver, obtaining solutions on certain irregular domains demands tailored approaches that are challenging to design and solve.

In recent year, deep learning approaches have become a promising and popular methodology for the numerical solution of various PDEs. They have the potential to overcome some of the challenges that classical methods are facing. That is, neural networks have the advantage to not generally rely on a grid. By leveraging automatic differentiation \cite{Baydin2018}, they eliminate the need for discretisation. Additionally, neural networks are able to represent more general functions than, e.g., an FEM basis, and they do not suffer from the curse of dimensionality. While the training of a neural network can become computationally demanding, especially when it consists of a non-convex optimisation problem, it is very efficient when evaluating new data points once trained. In this emerging field of deep learning for approximating PDE solutions, the class of approaches closest to the classical methodologies is the one of function approximators \cite{Maass2022}. They essentially model the PDE solution by a deep neural network and train the network's parameters to approximate the solution. Such approaches are, for example, the Deep Ritz method \cite{E2018} or the Deep Galerkin method \cite{Sirignano2018}. A widely used and adapted method of this class are the so-called \emph{physics-informed neural networks} (PINNs) \cite{Raissi2019} that will be the focus of this paper. Originally published in a two-part instalment \cite{Raissi2017a,Raissi2017b}, Raissi et al. developed the vanilla PINNs approach \cite{Raissi2019}. The basic idea behind PINNs is to minimise an energy functional that is the residual of the PDE and its initial and boundary conditions. The neural network itself models the solution function $u(t,x)$ given input variables $t$ and $x$ based on the underlying PDE. It has shown great success for many different types of PDEs \cite{PINNshighspeedflows,PINNsearthquakes} and has been extended to various specialised cases \cite{cPINNs,VPINNs,B-PINNs,surveyCuomo}.

While deep learning approaches for PDEs have gained a lot of traction in the last years and are being employed in increasingly more applications, they come with their set of challenges. In this work, we therefore systematically compare the finite element method and physics-informed neural networks in a computational study. Before giving a complete outline of the following, we review recent works on PINNs.

\subsection{PINNs: state of the art.}
As opposed to the finite element method, the theoretical groundwork for PINNs is rather sparse. The first work on convergence results with respect to the number of training points is by Shin et al. \cite{convEllipticParabolic}. For linear second-order elliptic and parabolic PDEs, they prove strong convergence in $C^0$ for i.i.d. sampled training data. Mishra et al. \cite{Mishra2022} have in turn developed upper bounds on the generalisation error of PINNs given some stability assumptions on the PDE. Focusing on a specific PDE, Ryck et al. \cite{errorNavierStokes} investigate the incompressible Navier-Stokes equation and provide an upper bound on the total error. However, their considerations are limited to the case of networks with two hidden layers and $\tanh$ activation functions. Despite the remaining need for more extensive theoretical work, PINNs have been extended into many different directions. Jagtap et al. \cite{cPINNs} develop conservative PINNs (cPINNs) to incorporate complex non-regular geometries by decomposing a spatial domain into independent parts and train separate PINNs for each. This work has been generalised to the extended PINNs (XPINNs) \cite{XPINNs} to allow space-time domain decomposition that can be applied to any type of PDE and enables parallelisation in training. Another domain decomposition method for PINNs are the hp-VPINNs \cite{hp-VPINNs}. This work is based on the variational PINNs \cite{VPINNs} that take inspiration from classical approaches for solving PDEs. VPINNs form the loss functional based on the variational form of the PDE with Legendre polynomials as test functions. Thus, allowing, e.g., for certain non-smoothnesses in PDEs and also for more efficient training. Bayesian PINNs \cite{B-PINNs} for noisy data employ a Bayesian neural network. Finite Basis Physics-Informed Neural Networks \cite{FBPINNs} combine the PINN-idea with the finite element method, aiming to reduce the spectral bias in PINNs \cite{spectralbias}. Interesting applications of PINNs appear in fluid dynamics \cite{PINNshighspeedflows}, electromagnetism \cite{PINNseigenvalueproblems}, elastic material deformation \cite{PINNselasticdeformation} and seismology \cite{PINNsearthquakes}. For a more extensive review of PINNs, its applications and extended forms, we refer to the survey by Cuomo et al. \cite{surveyCuomo}. Some shortcomings of PINNs are documented by Krishnapriyan et al. \cite{PINNsfailuremodes}, who find that PINNs struggle to learn relevant physics in more challenging PDE regimes. However, they simultaneously present ideas to address these issues.

PINNs were created with physical application in mind. Thanks to their flexibility, the main building block can be kept the same across many physical problems with smaller adjustments to the architecture. Implementation of the framework has also become more approachable since the release of dedicated software packages such as DeepXDE \cite{DeepXDE}, NVIDIA Modulus (previously SimNet) \cite{SimNet} or NeuroDiffEq \cite{NeuroDiffEq}. 

\subsection{Contributions and outline}
As mentioned above, the main goal of this work is a systematic comparison of physics-informed neural networks and the finite element method for the solution of partial differential equations. Indeed, we consider
\begin{itemize}
    \item the elliptic Poisson equation in one, two, and three space-dimensions,
    \item the parabolic Allen-Cahn equation in one space-dimension, and 
    \item the hyperbolic semilinear Schr\"odinger equation in one and two space dimensions.
\end{itemize}
We choose these model problems to cover a large range of classes of partial differential equations. We compare PINNs and FEM in terms of solution time, evaluation time, and accuracy. We employ different finite element bases, as well as a multitude of network architectures. Several improvements over the vanilla PINNs have been discussed in the literature. As those are either still fundamentally based on the same idea or rather represent a mix of a classical approach and PINNs, we focus on vanilla PINNs in the following.

This work is structured as follows. We discuss the PDEs, FEM, and PINNs in Section~\ref{sec_Math_back}. We introduce our method of comparison in Section~\ref{sec_method=comp}, before presenting our computational results regarding Poisson, Allen--Cahn, and semilinear Schr\"odinger equation in Sections~\ref{sec_Poisson}, \ref{sec_ACE}, and \ref{Sec_SSE}, respectively. We discuss our results and conclude the work in Section~\ref{sec_Discu}.

\section{Mathematical Background} \label{sec_Math_back}
In the following, we first study a very general class of partial differential equations that we formally denote by 
\begin{align} \label{eq.generalPDE}
    \mathcal{A}u(x,t) = f(x,t) \quad x \in \Omega, \,\, t \in [0,T].
\end{align}
Here, the function $u: \overline{\Omega} \times [0,T] \to \mathbb{R}^n$ denotes the solution of the PDE, where $\Omega \subset \mathbb{R}^d$ is open, bounded, and connected and usually represents a spatial domain, whereas $[0,T]$ is the time interval. $\mathcal{A}$ denotes the differential operator acting on $u$ that can also be nonlinear and $f: \Omega \times [0,T] \to \mathbb{R}^n$ is a source term. We denote the  corresponding boundary conditions and the initial condition by
\begin{align} \label{eq.generalBC}
     \mathcal{B}u(x,t) &= g(x,t), \quad x \in \partial \Omega, \,\, t \in [0,T] \\
    u(x,0) &= h(x), \quad x \in \Omega,
\end{align}
respectively.

Throughout this work, we consider three different partial differential equations: the \emph{Poisson equation}, the \emph{Allen--Cahn equation}, and the \emph{semilinear Schr\"odinger equation}. We introduce those PDEs in the following.

\paragraph{Poisson equation.} The Poisson equation is a linear elliptic PDE of the form
$$
\Delta u(x) = f(x), \qquad x \in \Omega,
$$
where $\Delta = \sum_{i=1}^d \frac{\partial}{\partial x_i} $ denotes the Laplacian and $f:\overline{\Omega} \rightarrow \mathbb{R}$ is a source term. To be well-defined, we need to equip it with boundary conditions, say, of \emph{Dirichlet}-type:
$$
u(x) = u_\partial(x), \qquad x \in \partial \Omega,
$$
\emph{Neumann}-type:
$$
\partial_{\vec{n}} u(x) = u_\partial(x) \qquad x \in \partial \Omega,
$$
or a combination of the two. In each case, we employ the function $u_\partial: \partial \Omega \rightarrow \mathbb{R}$ to determine the boundary behaviour. There are several results about the existence of strong and weak solutions of elliptic PDEs, see, e.g. Theorem 3 in Chapter 6.2 in \cite{Evans}. The Poisson equation models, for instance, the stationary heat distribution in a homogeneous object $\Omega$; here, $f$ describes heat sources and sinks.

\paragraph{Allen--Cahn equation.} The Allen--Cahn equation is a nonlinear parabolic PDE of the form
$$
\frac{\partial u(x,t)}{\partial t}= \varepsilon\Delta u(x,t) - \frac{2}{\varepsilon}u(t,x) (1-u(t,x)) (1-2u(t,x)), \qquad x \in \Omega, t \in [0, T],
$$
where $\varepsilon>0$ and which, of course, is considered together with appropriate boundary conditions and an initial condition. The PDE is semilinear, since the nonlinearity depends on $u$ but not on a second derivative of $u$. If $\varepsilon$ is sufficiently small, the solution of the Allen--Cahn equation approximately partitions the domain $\Omega$ into patches where $u \approx 0$ and $u\approx 1$; in-between those patches -- at the so-called diffuse interface -- it is smooth. These patches can represent binary alloys \cite{Allen-Cahn} or, e.g., different segments in an image with pixels in $\Omega$ \cite{Benes}. Allen--Cahn is also a relaxation of the mean-curvature flow to which it converges as $\varepsilon \downarrow 0$ \cite{Feng}.

\paragraph{Semilinear Schr\"odinger equation.} The semilinear Schr\"odinger equation is a complex-valued nonlinear hyperbolic PDE of the form
$$
\mathrm{i}\frac{\partial h(x,t)}{\partial t} = - \Delta h(x,t) - h(x,t)|h(x,t)|^2, \qquad x \in \Omega, t \in [0, T],
$$
again subject to appropriate initial and boundary conditions. We have represented the semilinear Schr\"odinger equation above as it is usual in the literature. This way of writing might be confusing. Thus, we present the PDE again splitting real and imaginary parts. We set  $h(x,t) =: u_R(x,t) + \mathrm{i}u_I(x,t)$ and write
\begin{align*}
 \frac{\partial u_R(x,t)}{\partial t} &= - \Delta u_I(x,t) - (u_R(x,t)^2 +u_I(x,t)^2)u_I(x,t),\\
 \frac{\partial u_I(x,t)}{\partial t} &= \Delta u_R(x,t) + (u_R(x,t)^2 +u_I(x,t)^2)u_R(x,t), \qquad x \in \Omega, t \in [0, T].
\end{align*}
We refer to \cite{STRAUSS} for an application of the semilinear Schr\"odinger equation in non-linear optics.

\subsection{Finite Element Method}\label{sec.FEM}
As mentioned before, the finite element method has been the gold standard for the spatial discretisation of a huge class of partial differential equations. We now discuss the basics of the finite element method given the example of an elliptic PDE with Dirichlet boundaries $u_\partial = 0$ and square intregable $f \in \mathcal{L}^2(\Omega)$. We mention non-stationary PDEs and other boundary conditions further below.

When solving an elliptic PDE with the finite element method, we  aim to find a weak solution. That means, we try to find $u$ in an appropriate function space $U$ such that for all functions $v$ in an appropriate $V$, we have
$$
\int_\Omega v(x) (\Delta u(x) -  f(x)) \mathrm{d}x = 0.
$$
Applying integration by parts, we obtain the usual weak formulation of the Poisson equation:

\begin{equation} \label{eq:var:poisson}
    \int_\Omega  \langle \nabla v(x), \nabla u(x)\rangle +  v(x)f(x) \mathrm{d}x = 0.
\end{equation}

An appropriate choice of function spaces $U,V$ in this case is $U = V = H_0^1(\Omega)$, the Sobolev space of square-integrable functions $\Omega \rightarrow \mathbb{R}$ that have a weak derivative that is also square-integrable. 
In finite elements, we now replace $H_0^1(\Omega)$ by a finite-dimensional space $H$ with basis $(\varphi_i)_{i=1}^N$. On this finite dimensional space, we can write the weak equation as 

$$
-\int_\Omega  \langle \nabla \varphi_i(x), \nabla \sum_{j=1}^n a_j\varphi_j(x)\rangle \mathrm{d}x =  \int_\Omega \varphi_i(x)f(x) \mathrm{d}x \qquad i=1,...,N,
$$
which can be rewritten as a usual linear system  $Ba = b$, with
$$
B = \left(-\int_\Omega  \langle \nabla \varphi_i(x), \nabla \varphi_j(x)\rangle \mathrm{d}x \right)_{i,j=1}^N, \qquad b = \left(\int_\Omega \varphi_i(x)f(x) \mathrm{d}x \right)_{i=1}^N.
$$
Now, of course, $u(\cdot) \approx \sum_{j=1}^n a_j\varphi_j(\cdot)$. While this approach allows for a large range of bases $(\varphi_i)_{i=1}^N$, we usually speak only of finite elements when choosing specific locally supported, piecewise polynomial functions. As usual differential operators are local, this will usually lead to $B$ having a favourable, sparse structure.

In the discussion above, we consider homogeneous Dirichlet boundary conditions $u_\partial = 0$. Inhomogeneous boundary conditions, in which $u_\partial$ is not constantly $0$, can be achieved by first finding a function that satisfies the boundary conditions and then solving an auxiliary homogeneous problem with $u_\partial = 0$. Neumann conditions can be enforced by adding to the bilinear form $\int_{\partial \Omega} v u_\partial \mathrm{d}x$ and choosing functions in $U = V = H^1(\Omega)$ the space of square-integrable functions with square-integrable weak derivatives. Mixed boundary conditions can be enforced in a similar way.

The time-domain of a non-stationary PDE can also be discretised with finite elements, see, e.g., \cite{Hulbert}. However, the probably more usual approach is to use finite elements in space and usual ODE solvers in time. To represent Allen--Cahn and semilinear Schr\"odinger at the same time, we consider some semilinear PDE of the form
$$
\frac{\partial u(x,t)}{\partial t} = \Delta u(x,t) + F(u(x,t)), \quad u(x,0) = u_0(x) \qquad (x \in \Omega, t\in[0,T]),
$$
which is also subject to boundary conditions. Due to the stiffness of the heat equation, we are required to use implicit techniques, such as the implicit Euler method. In a semi-discrete form, we can write the update step as
$$u_{t+1} = u_t + \delta t \Delta u_{t+1} + \delta t F(u_{t+1}) \qquad (t = 0, 1, \ldots),$$
where $\delta t>0$ denotes the  size of the time steps.
In the weak formulation, we obtain 
$$\int_{\Omega} v u_{t+1} \mathrm{d}x = \int_{\Omega} v u_t - \delta t\langle \nabla v, \nabla u_{t+1}\rangle + \delta t v F(u_{t+1}) \mathrm{d}x \qquad (v \in V, t = 0, 1, \ldots)$$
and can now again use a finite element discretisation in space to obtain an approximation of $u_{t+1}$ represented through a finite-dimensional equation. In the case of Allen--Cahn and semilinear Schrödinger, this equation is non-linear and requires us to repeatedly employ a Newton's method. 

We can prevent the additional cost of a Newton-solve and usually still obtain a stable discretisation, by employing a semi-implicit strategy: linear parts of the PDE are solved with an implicit discretisation, non-linear parts are discretised explicitly. In our setting, we obtain
$$u_{t+1} = u_t + \delta t \Delta u_{t+1} + \delta t F(u_{t}) \qquad (t = 0, 1, \ldots)$$
which requires us to solve the following weak problem 
$$\int_{\Omega} v u_{t+1}+ \delta t  \langle \nabla v, \nabla u_{t+1}\rangle \mathrm{d}x  = \int_{\Omega} v u_t + \delta tv F(u_{t}) \mathrm{d}x \qquad (v \in V, t = 0, 1, \ldots)$$
that is fully linear.

For more details on the finite element method, we refer to, e.g., the book by Braess \cite{braess_2007} or the classical references Courant \cite{Courant} and Hrennikoff \cite{Hrennikoff}. ODE integrators, as well as evolution equations, are thoroughly considered by Iserles \cite{iserles_2008}. Semi-implicit schemes appear, for instance, in the work by Bertozzi and Sch\"onlieb \cite{bertozzi_schonlieb_2010}.

\subsection{Physics-informed Neural Networks}\label{sec.PINNs}
The aim of physics-informed neural networks is to approximate PDE solutions using a deep neural network. They make use of the powerful tool that is automatic differentiation and therefore do not rely on discretisation of the space-time domain but rather on random sampling of the domain. 

Let us consider a PDE of the form \eqref{eq.generalPDE} with suitable boundary and initial conditions \eqref{eq.generalBC}. The vanilla PINNs approach, as introduced by Raissi et al. \cite{Raissi2019}, uses a fully connected neural network with $N_l$ hidden layers and $N_e(l)$ neurons per layer $l$. Inputs of the network are the PDE variables $(x,t)$ sampled in the domain $\Omega \times [0,T]$. The neural network acts as the function approximator $u_{\theta}(x,t)$ of the PDE solution, with $\theta$ the network weights that are optimised during training. The PDE is integrated as a soft constraint in the optimisation. That is, the network is trained with the PDE residual, as well as boundary and initial condition residuals as the loss functional:
\begin{align*}
    Loss(\theta) = &\frac{1}{N_f} \sum_{i = 1}^{N_f} \lVert \mathcal{A}u_{\theta}(x^i_f,t^i_f) - f(x^i_f,t^i_f)) \rVert^2 \\
    &+ \frac{1}{N_g} \sum_{j = 1}^{N_g} \lVert \mathcal{B}u_{\theta}(x^j_g,t^j_g) - g(x^j_g,t^j_g) \rVert^2 + \frac{1}{N_h} \sum_{k = 1}^{N_h} \lVert u_{\theta}(x^k_h,0) - h(x^k_h) \rVert^2.
\end{align*}
Here, $N_f$ is the number of collocation points $(x_f^i, t_f^i) \in \Omega \times [0,T]$ for $i = 1, \dots, N_f$ sampled for the PDE residual in the loss. Similarly, $(x_g^j,t_g^j) \in \partial \Omega \times [0,T]$ for $j = 1, \dots, N_g$ denote the training points on the boundary and $x_k^h \in \Omega$ for $k = 1, \dots, N_h$ the training data for the initial condition. Additionally, we need data $g(x^j_g,t^j_g)$ for the boundary and $h(x^k_h)$ initial conditions. This can be measured data and does not need to be represented as an analytic function. PINNs is therefore able to incorporate measurements and leverage data-driven information. We follow \cite{Raissi2019} in using Latin Hybercube Sampling \cite{Stein1987}, a quasi-random approach for space filling sampling, to obtain the collocation points for training. In our experiments, we re-sample the collocation points in every epoch to get a better coverage of the sampling domain; see also \cite{Jin}. As mentioned above, the differential operator $\mathcal{A}$ and any derivatives in the boundary condition are evaluated using automatic differentiation \cite{Baydin2018}. In contrast to numerical methods such as finite differences, automatic differentiation uses the chain rule to backpropagate through the network and evaluate the derivative. It is therefore not dependent on a grid or mesh that discretises the domain. In turn, the sampling method becomes more important. Automatic differentiation gives rise to a significant advantage of PINNs towards other classical numerical methods for solving PDEs, that is, it does not deal with prohibitively small step-sizes and scales well in higher dimensions. In PINNs, the design of the neural network architecture, i.e., the type of network structure, the number of hidden layers and the number of nodes per layer, is flexible and can be adjusted based on the PDE complexity. In the following, we will use fully connected feed-forward dense neural network. We denote the size and numbers of nodes per layers as $[N_e(1),N_e(2),...,N_e(l)]$. That is, $[5,5,1]$, for example, represents a neural network with 2 hidden layers and 5 nodes per layer. The last entry 1 represents the size of the output layer. Lastly, the loss function is typically minimised using first the Adam optimiser \cite{Kingma2015} for a coarser optimisation and then the second-order quasi-Newton optimiser L-BFGS \cite{Liu1989}. 

\section{Method of comparison} \label{sec_method=comp}
In this section, we give an overview of the experimental design and the measures that we consider to compare FEM and PINNs. The main features that we investigate are the time to solve the PDEs and the accuracy of the results. We therefore ask the questions of which methodology is computationally faster, more accurate, and most efficient. We investigate different types of PDEs with varying complexity and dimensionality to cover a broad spectrum of PDEs and examine if computational speed and accuracy of the two approaches change based on the PDE type.

\subsection{Experimental Setup}
We now discuss the experimental setup of our computational study. We first discuss ground truth solutions, and then describe the setup for FEM and PINNs and the metrics with which compare them. 

Indeed, for a systematic comparison of FEM and PINNs, we need a ground truth solution of the PDEs to compare the solution approximations to and evaluate their accuracy. The first step is therefore to determine these ground truth solutions. For the Poisson equations, we have analytical solutions that we can use to determine the accuracy of the approximation approaches. However, neither the Allen-Cahn equation nor the semilinear Schrödinger equation have analytical solutions. Instead, we use FEM on a very fine mesh to compute a very accurate reference solution that we will use as the ground truth; time-stepping is performed using the implicit Euler method. One of the advantages of FEM is that we have extensive theoretical foundations for the convergence for large classes of partial differential equations, including the Allen--Cahn equation, e.g., \cite{Feng2005} and the semilinear Schr\"odinger equation, e.g. \cite{Shi}. Thus, a finely meshed finite element solution can therefore guarantee accuracy up to a small error.

Let us now discuss the details and specifications for FEM that we consider in the comparison. As mentioned before, finer meshes will lead to higher accuracy albeit slower computation time. We therefore solve the PDEs for different mesh sizes to see the relationship between computation time and accuracy based on the FEM design. For time-dependent PDEs, that are, the Allen-Cahn and the semilinear Schrödinger equations, we choose a semi-implicit strategy for discretisation as detailed in Subsection \ref{sec.FEM}. All FEM solution approximations are implemented using the python toolbox FEniCS \cite{Alnaes2015,Logg2012}. In the case of PINNs, we closely follow the vanilla approach as described in Raissi et al. \cite{Raissi2019}. That is, we design the loss functional as in Subsection \ref{sec.PINNs} and use the Adam optimiser \cite{Kingma2015} in the first instance and L-BFGS \cite{Liu1989} to refine the optimisation. The learning rate is heuristically chosen for optimal results in each PDE. All derivatives in the loss functional are computed using automatic differentiation. The collocation points that make up the input to the neural network are newly sampled for each epoch using Latin Hybercube sampling \cite{Stein1987} in the Adam optimisation. This way, we are able to cover more of the sampling space and improve generalisability of the trained network. L-BFGS does not allow for re-sampling between iterations and we therefore only sample the collocation points once before the optimisation again using Latin Hybercube sampling. PINNs allow for a flexible design of the neural network architecture based on the underlying PDE. We therefore run the network training on architectures of different sizes, that is, with varying numbers of layers and nodes. The code for PINN training and evaluation was written with the python neural network library jax \cite{Jax2018} and is available on github: \url{https://github.com/TamaraGrossmann/FEM-vs-PINNs}. 

We consider three main features for comparison. That is, we evaluate the FEM and PINNs approaches based on their solution time, evaluation time, and accuracy. The solution time refers to the time it takes for each method to approximate the general solution. We differentiate between solution and evaluation time due to the way FEM and PINNs approximate solutions differently. FEM, typically, solve the PDE on a fixed mesh of interest but it is possible to interpolate that solution to a different mesh. On the other hand, for PINNs a neural network is trained on a set of collocation points. The trained network can then be evaluated on any point in the domain. While training time of PINNs can be rather slow, one of the advantages of neural networks is that once trained the evaluation on new input data tends to be very fast. We therefore consider both times. For FEM this means solution time refers to solving the weak form of the PDE on a fixed mesh. FEM is run on a CPU. Considering PINNs, the solution time refers to the training time of the neural network which is run on a GPU. In turn, the evaluation time in FEM is measured as the time to interpolate the solution on a different mesh. In PINNs, the evaluation time is the time to evaluate the trained neural network on a new set of collocation points. We measure the accuracy of both methods on the same mesh in comparison to the ground truth solutions. For the Allen-Cahn and Schrödinger equations the evaluation mesh is chosen as the fine mesh on which the ground truth solution was derived. The measure for accuracy is the $\ell^2$ relative error. All experiments were run on the same machine that has 12 CPU cores and we used a Quadro P6000 GPU for the neural network training. The codes were run 10 times and the reported solution and evaluation times are the average over the 10 recorded times.

\section{Approximating the Poisson equation} \label{sec_Poisson}
Let us first investigate the Poisson equation in one, two, and three space-dimensions. Each of the equations we consider has an analytical solution that we can use to evaluate the accuracy of the FEM and PINN approximation. Comparing the same type of PDE in different dimensions also allows us to draw conclusions about cost and accuracy effects due to the dimensionality of the PDE.

\subsection{1D}
\begin{figure*}[t]
    \centering
     \begin{subfigure}[b]{0.49\textwidth}
         \centering
        \includegraphics[width=\textwidth]{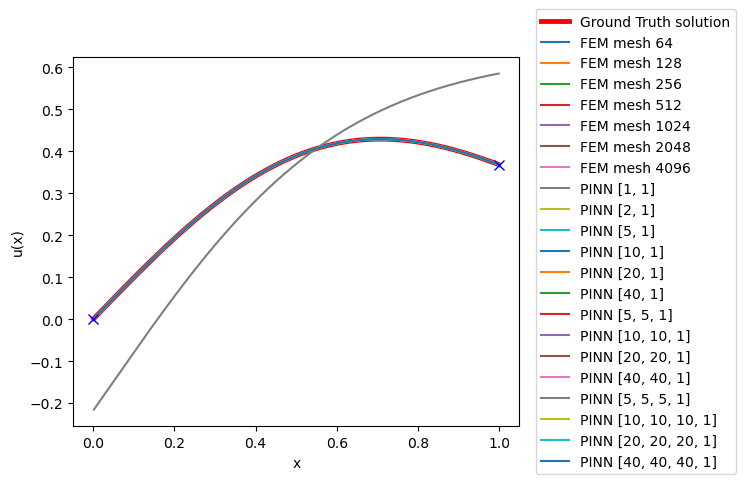}
        \caption{Plot of the PDE ground truth solution and the FEM and PINNs solution approximations.}
         \label{fig.1DPoisson_solution}
     \end{subfigure}
     \hfill
     \begin{subfigure}[b]{0.49\textwidth}
         \centering
         \includegraphics[width=\textwidth]{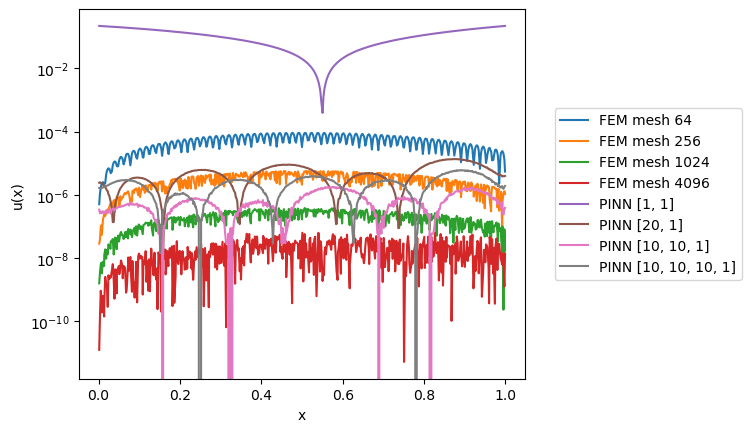}
         \caption{Plot of the difference between some of the FEM / PINNs solution approximations to the ground truth solution on a log scale.}
         \label{fig.1DPoisson_solution_difference}
     \end{subfigure}
        \caption{Plot for 1D Poisson equation solution.}
\end{figure*}

For the one-dimensional case, we are examining the Poisson equation with a right hand side $f(x)$ as detailed below on a unit interval and employ Dirichlet boundary conditions:
\begin{align} \label{eq.1DPoisson}
    \begin{split}
    \Delta u(x) &= (4x^3 - 6x) \exp{(-x^2)} \quad x \in (0,1) \\
    u(0) &= 0 \\
    u(1) &= \exp(-1). 
    \end{split}
\end{align}
This PDE has an analytical solution that can be written as:
\begin{align*}
    u_{\rm true}(x) = x\exp{(-x^2)} \qquad (x \in [0,1]).
\end{align*}
A visualisation of the solution to the 1D Poisson equation is shown in Figure \ref{fig.1DPoisson_solution}.

\subsubsection*{FEM}
As introduced in Section \ref{sec.FEM}, the first step in solving PDEs with the finite element method is deriving the weak formulation of the PDE, which we have done for Poisson equation alright right there. Here, of course $f(x) = (4x^3 - 6x) \exp{(-x^2)}$. Next, we need to define the finite element mesh. We choose a regular mesh on the domain $[0,1]$ with varying numbers of cells $n \in \{64,128,256,512,1024,2048,4096\}$. Of course, a higher number of cells corresponds to a finer grid on which the PDE is solved and subsequently leads to a more accurate, but also  computationally more costly solution.  The finite elements we are using are standard linear Lagrange elements, that is,  piecewise linear hat functions ($P_1$). Subsequently, we solve the variational problem in \eqref{eq:var:poisson} with the Dirichlet boundary conditions specified in \eqref{eq.1DPoisson} using the conjugate gradient method with incomplete LU factorisation as a preconditioner. All specifications are implemented using the python toolbox FEniCS \cite{Alnaes2015,Logg2012} to compute an approximate PDE solution.

\subsubsection*{PINNs}
For solving the 1D Poisson equation using PINNs, there are three design parameters that we need to specify before training. The first step is choosing a loss functional. Following the vanilla PINNs approach, we evaluation the goodness of the solution using the discretised $\ell^2$-energy or mean squared error over the PDE, boundary and initial conditions. In particular, we define the loss function as:
\begin{align*}
    \text{Loss}(\theta) = &\frac{1}{N} \sum_{i = 1}^{N} \lVert 
    \Delta u_{\theta}(x_i) - (4x_i^3 - 6x_i) \exp{(-x_i^2)} \rVert^2_2 
    + \lVert u_{\theta}(0) \rVert^2_2  + \lVert u_{\theta}(1) - exp(-1) \rVert^2_2,
\end{align*}
with $u_{\theta}$ the neural network, $\theta$ the trained weights and $N = 256$ the number of collocation points $x_i$ sampled in each epoch using latin hybercube sampling. The second design parameter is the neural network architecture, that is, the type of neural network, the activation function, and the number of hidden layers and nodes. For the 1D Poisson case, we train feed-forward dense neural networks with $\tanh$ as the activation function. We compare the results on architectures of different sizes. The network architectures we consider for 1D Poisson are $[1,1], [2,1], [5,1], [10,1], [20,1], [40,1], [5,5,1], [10,10,1],$ $[20,20,1], [40,40,1], [5,5,5,1], [10,10,10,1], [20,20,20,1],$ and $[40,40,40,1]$. The loss function is minimised for each network architecture using the Adam optimiser for $15,000$ epochs with a learning rate of $1e-4$ in the first instance. Additionally, we refine the optimisation using L-BFGS.

\subsubsection*{Results}
\begin{figure*}[t]
     \centering
     \begin{subfigure}[b]{0.49\textwidth}
         \centering
         \includegraphics[width=\textwidth]{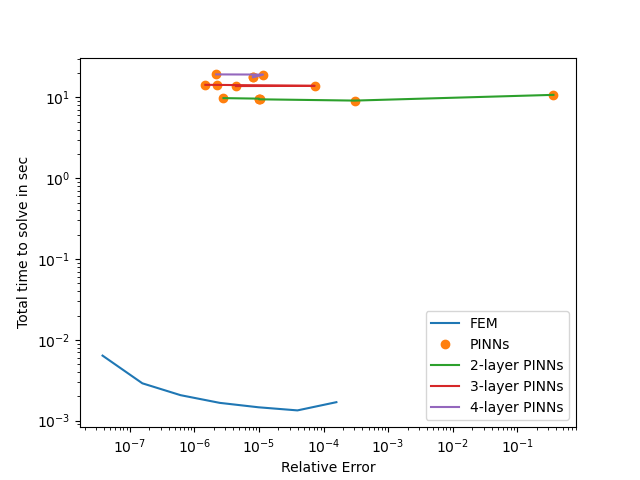}
         \caption{Plot of time to solve FEM and train PINN in sec versus $\ell^2$ relative error.}
         \label{fig.1DPoisson_solve}
     \end{subfigure}
     \hfill
     \begin{subfigure}[b]{0.49\textwidth}
         \centering
         \includegraphics[width=\textwidth]{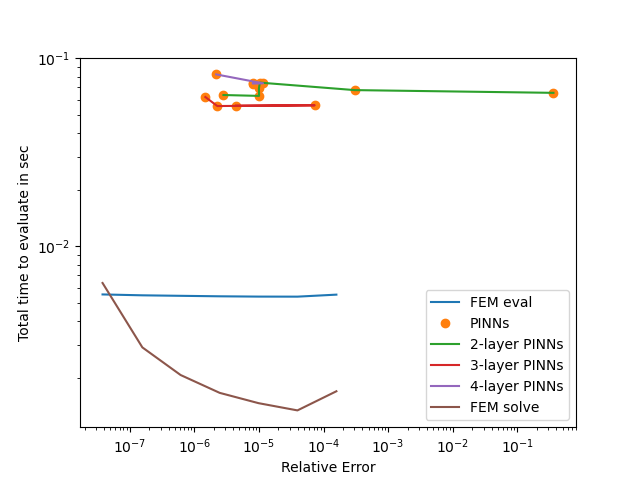}
         \caption{Plot of time to interpolate FEM and evaluate PINN in sec versus $\ell^2$ relative error.}
         \label{fig.1DPoisson_test}
     \end{subfigure}
        \caption{Plot for 1D Poisson equation of time in sec versus $\ell^2$ relative error.}
        \label{fig.1DPoisson_plots}
\end{figure*}
The resulting PDE solution approximations of the 1D Poisson equation using FEM and PINNs are compared to the analytical solution on a mesh in $[0,1]$ with 512 mesh points. The ground truth solution of the 1D Poisson equation in \eqref{eq.1DPoisson} and its approximations are displayed in Figure \ref{fig.1DPoisson_solution}. Likewise, the difference from the approximate solutions to the GT solution are shown in Figure \ref{fig.1DPoisson_solution_difference} for all mesh sizes of FEM and architectures in PINNs. One of the architectures in the PINNs approximation renders results with large relative error: the PINN with a single hidden layer and one node is not even able to learn the solution in a way in which the boundary conditions are satisfied. However, all other solution approximations differ from the ground truth only marginally.

Let us compare the time it takes to solve or rather approximate the PDE using FEM and PINNs to the relative error produced on a new set of grid points. For FEM the solution time is the time to solve the linear systems and for PINNs we consider the time to train the neural network. The results are shown in Figure \ref{fig.1DPoisson_solve}. We can clearly show that overall FEM are faster and more accurate in their solution approximation. While there are some PINN architectures that are able to achieve similar or even lower relative errors than some of the FEM approximations, their training time is 2-3 orders of magnitude higher than in FEM. Considering the evaluation time, that is, the time to interpolate the FEM solution on a different mesh and evaluate the trained PINN on a test set, we can see a similar relationship; see Figure \ref{fig.1DPoisson_test}. FEM solution approximations are overall faster and more accurate. Finally, we consider the relationship between the number of layers and the solution time and relative error. We observe that the time to train a PINN is similar relative to the number of layers. However, the accuracy of each network is related to the number of nodes in the layers. That is, we cannot show that networks with more layers generally achieve better results in 1D Poisson.

\subsection{2D}
Let us now investigate the two-dimensional Poisson equations defined as:
\begin{align}\label{eq.2DPoisson}
    \begin{split}
    \Delta u(x,y) = 2& (x^4 (3 y - 2) + x^3 (4 - 6 y) + x^2 (6 y^3 - 12 y^2 + 9 y - 2) \\ 
    & - 6 x (y - 1)^2 y + (y - 1)^2 y) \qquad \qquad \qquad \qquad (x,y) \in (0,1)^2 
    \end{split} \\
    \begin{split} \label{eq.2DPoisson_bn}
   \partial_{\vec{n}} u(0,y) &= 0 \quad y \in [0,1]\\
    \partial_{\vec{n}} u(1,y) &= 0 \quad y \in [0,1]\\
    u(x,0) &= 0 \quad x \in [0,1]\\
   \partial_{\vec{n}}  u(x,1) &= 0 \quad x \in [0,1]
   \end{split}
\end{align}
We can solve the PDE \eqref{eq.2DPoisson} with mixed boundary conditions analytically:
\begin{align*}
    u_{\rm true}(x,y) = x^2  (x-1)^2 y (y-1)^2.
\end{align*}
The solution is displayed in Figure \ref{fig.2DPoisson_solution}.

\subsubsection*{FEM}
For the 2D Poisson equation, we again refer to the weak formulation of the Poisson equation in \eqref{eq:var:poisson} with $f(x,y) = 2 (x^4 (3 y - 2) + x^3 (4 - 6 y) + x^2 (6 y^3 - 12 y^2 + 9 y - 2)  - 6 x (y - 1)^2 y + (y - 1)^2 y)$. The finite element mesh is defined on the unit square $[0,1] \times [0,1]$; it contains $n \times n$ squares wih $n \in \{100, 200,\ldots,1000\}$. Each square on the mesh is divided into two triangles on which we define again piecewise linear finite elements ($P_1$). We solve the variational problem \eqref{eq:var:poisson} with mixed boundary conditions given in \eqref{eq.2DPoisson_bn} using the conjugate gradient method with incomplete LU factorisation preconditioning. The method is implemented using FEniCS.
 
\subsubsection*{PINNs}
The loss functional for the PINN approximation is the $\ell^2$-residual of the PDE \eqref{eq.2DPoisson} and its boundary conditions. For $u_{\theta}$ the neural network with weights $\theta$ that are to be trained, the loss reads:
\begin{align*}
    \text{Loss}&(\theta) \\ = \, &\frac{1}{N_f} \sum_{i = 1}^{N_f} \lVert 
     \Delta u_{\theta}(x^i_f,y^i_f) - 2 ((x^i_f)^4 (3 y^i_f - 2) + (x^i_f)^3 (4 - 6 y^i_f) + (x^i_f)^2 (6 (y^i_f)^3 - 12 (y^i_f)^2 + 9 y^i_f - 2) \rVert^2_2 \\
    &+ \frac{1}{N_g} \sum_{j = 1}^{N_g}\left(\left\lVert \partial_{\vec{n}}u_{\theta} (0,y^j_g)\right\rVert^2_2 + \left\lVert \partial_{\vec{n}} u_{\theta} (1,y^j_g)\right\rVert^2_2 + \left\lVert u_{\theta}(x^j_g,0)\right\rVert^2_2  + \left\lVert \partial_{\vec{n}} u_{\theta}(x^j_g,1)\right\rVert^2_2 \right).
\end{align*}
The collocation points are sampled at every epoch using Latin Hybercube sampling with $N_f = 2000$ and $N_g = 250$. We train a feed-forward dense neural network with $\tanh$ activation function. We consider 11 different architectures for training, these are [20,1], [60,1], [20,20,1], [60,60,1], [20,20,20,1], [60,60,60,1], [20,20,20,20,1], [60,60,60,60,1], [20,20,20,20,20,1], [60,60,60,60,60,1], and [120,120,120,120,120,1]. Just like in 1D Poisson, we use the Adam optimiser for $20,000$ epochs with a learning rate of $1e-3$ to train the network. Subsequently, we refine the optimisation using L-BFGS.

\begin{figure*}[t]
    \centering
    \includegraphics[width=\textwidth]{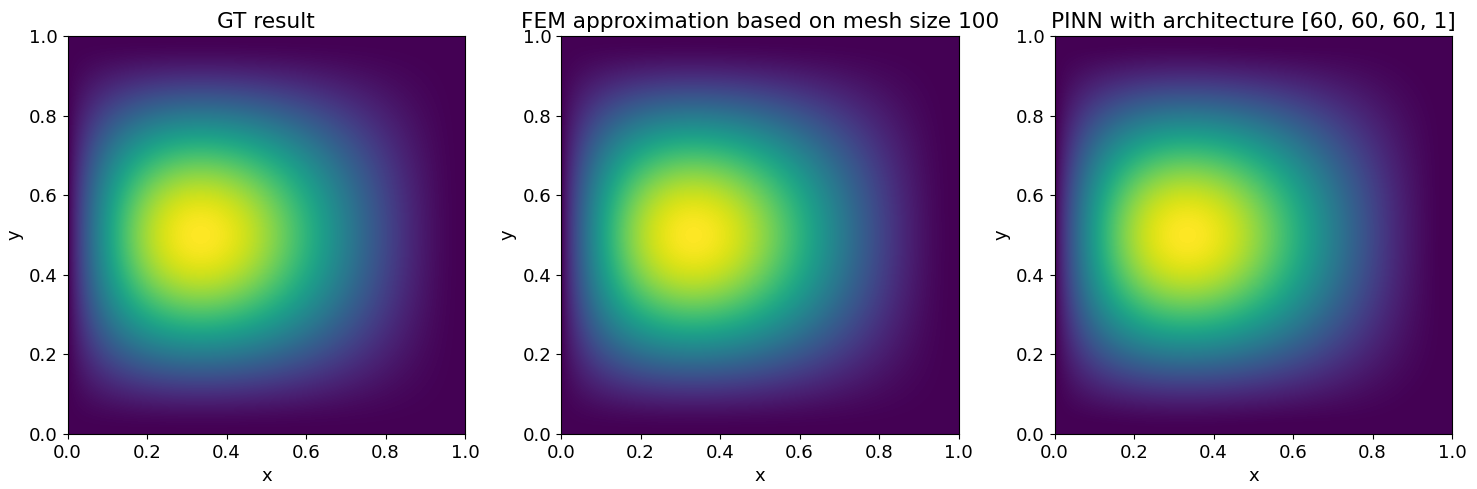}
    \caption{Comparison of the 2D Poisson ground truth solution to examples of the FEM and PINN approximations.}
    \label{fig.2DPoisson_solution}
\end{figure*}

\subsubsection*{Results}
\begin{figure*}[t]
     \centering
     \begin{subfigure}[b]{0.43\textwidth}
         \centering
         \includegraphics[width=\textwidth]{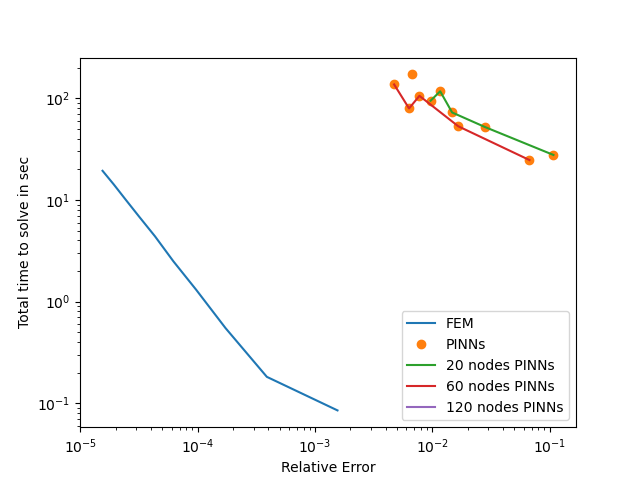}
         \caption{Plot of time to solve FEM and train PINN in sec versus $\ell^2$ relative error.}
         \label{fig.2DPoisson_solve}
     \end{subfigure}
     \hfill
     \begin{subfigure}[b]{0.43\textwidth}
         \centering
         \includegraphics[width=\textwidth]{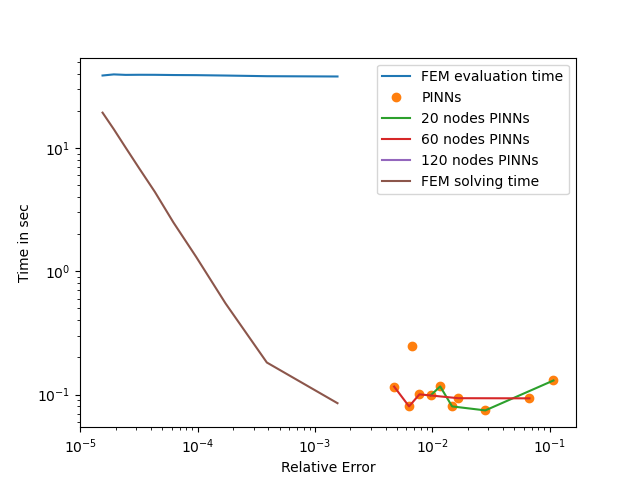}
         \caption{Plot of time to interpolate FEM and evaluate PINN in sec versus relative error. For comparison, the time to solve FEM is also plotted.}
         \label{fig.2DPoisson_test}
     \end{subfigure}
        \caption{Plot for 2D Poisson equation of time in sec versus $\ell^2$ relative error.}
        \label{fig.2DPoisson_plots}
\end{figure*}
An example of the resulting solution approximations on a mesh with $2000 \times 2000$ cells next to the ground truth solution is show in Figure \ref{fig.2DPoisson_solution}. The time versus accuracy plots for the 2D Poisson equation are displayed in Figure \ref{fig.2DPoisson_solve} and Figure \ref{fig.2DPoisson_test}. Considering the solution time, FEM clearly outperforms all PINN approximations both in accuracy and computation time. Using FEM to obtain the PDE solution is faster by 1-3 orders of magnitude. However, when we look at the evaluation times in Figure \ref{fig.2DPoisson_test}, the relationship changes. That is, the time to evaluate a PINN is 2-3 orders of magnitude faster then the interpolation of FEM on a new mesh. We additionally plotted the FEM solution time in the same plot, and the PINN evaluation remains faster, however, by a smaller margin. Interestingly, the solution time of FEM is faster than the FEM evaluation time; this is possibly due to an inefficient interpolation code. Even though the evaluation of the trained neural network gives an improvement in time, the PINN approximations remain to have lower accuracy.

\subsection{3D}
For the three-dimensional Poisson equation, we choose a PDE on the unit cube with Dirichlet boundary conditions, as follows:
\begin{align*}
    \begin{split}
   \Delta u(x,y,z) &=  -3 \pi^2\sin(\pi x)\sin(\pi y)\sin(\pi z) \quad (x,y,z) \in (0,1)^3 \\
    u(x,y,z) &= 0 \quad (x,y,z) \in  \partial(0,1)^3
    \end{split}
\end{align*}
The analytical solution of this 3D Poisson equation is displayed in Figure \ref{fig.3DPoisson_solution} and can be written as:
\begin{align*}
    u_{\rm true}(x,y,z) = \sin(\pi x)\sin(\pi y)\sin(\pi  z).
\end{align*}

\subsubsection*{FEM}
The weak formulation of the 3D Poisson equation is again given in \eqref{eq:var:poisson} with $f(x,y,z) = -3 \pi^2\sin(\pi x)\sin(\pi y)\sin(\pi z)$. The finite element mesh is defined on a unit cube $[0,1] \times [0,1] \times [0,1]$ consisting of $n \times n \times n$ cubes with $n \in \{16, 32, 64, 128\}$. We subdivide each cube into tetrahedrals and use again piecewise linear finite elements ($P_1$). The weak problem  is solved  using the conjugate gradient method and incomplete LU factorisation as the preconditioner.
 
\subsubsection*{PINNs}
Similar to the one- and two-dimensional Poission cases, we design the loss functional as the PDE and boundary condition residual. That is, we define the loss as follows:
\begin{align*}
    \text{Loss}(\theta) = \, &\frac{1}{N_f} \sum_{i = 1}^{N_f} \lVert 
     \Delta u_{\theta}(x^i_f,y^i_f,z^i_f) + 3 \pi^2\sin(\pi x^i_f)\sin(\pi y^i_f)\sin(\pi z^i_f) \rVert^2_2 \\
    &+ \frac{1}{N_g} \sum_{j = 1}^{N_g}\left(\lVert u_{\theta} (0,y^j_g,z^j_g)\rVert^2_2 + \lVert u_{\theta} (x^j_g,0,z^j_g)\rVert^2_2 + \lVert u_{\theta} (x^j_g,y^j_g,0)\rVert^2_2\right)\\
    &+\frac{1}{N_g} \sum_{j = 1}^{N_g}\left(\lVert u_{\theta} (1,y^j_g,z^j_g)\rVert^2_2 + \lVert u_{\theta} (x^j_g,1,z^j_g)\rVert^2_2 + \lVert u_{\theta} (x^j_g,y^j_g,1)\rVert^2_2\right).
\end{align*}
Here, $u_{\theta}$ denotes the neural network with $\theta$ the training parameters. We sample $N_f = 1000$ collocation points in the domain and $N_g = 100$ on the boundary. The neural network is again a feed-forward dense neural network with $\tanh$ activation function. We consider 8 different architectures to train with the following layer and node specifications: [20,20,1], [60,60,1], [20,20,20,1], [60,60,60,1], [20,20,20,20,1], [60,60,60,60,1], [20,20,20,20,20,1], and [60,60,60,60,60,1]. The Adam optimiser with learning rate $1e-3$ is used for 20,000 epochs before employing the L-BFGS optimiser.

\begin{figure*}[t]
    \centering
    \includegraphics[width=0.7\textwidth]{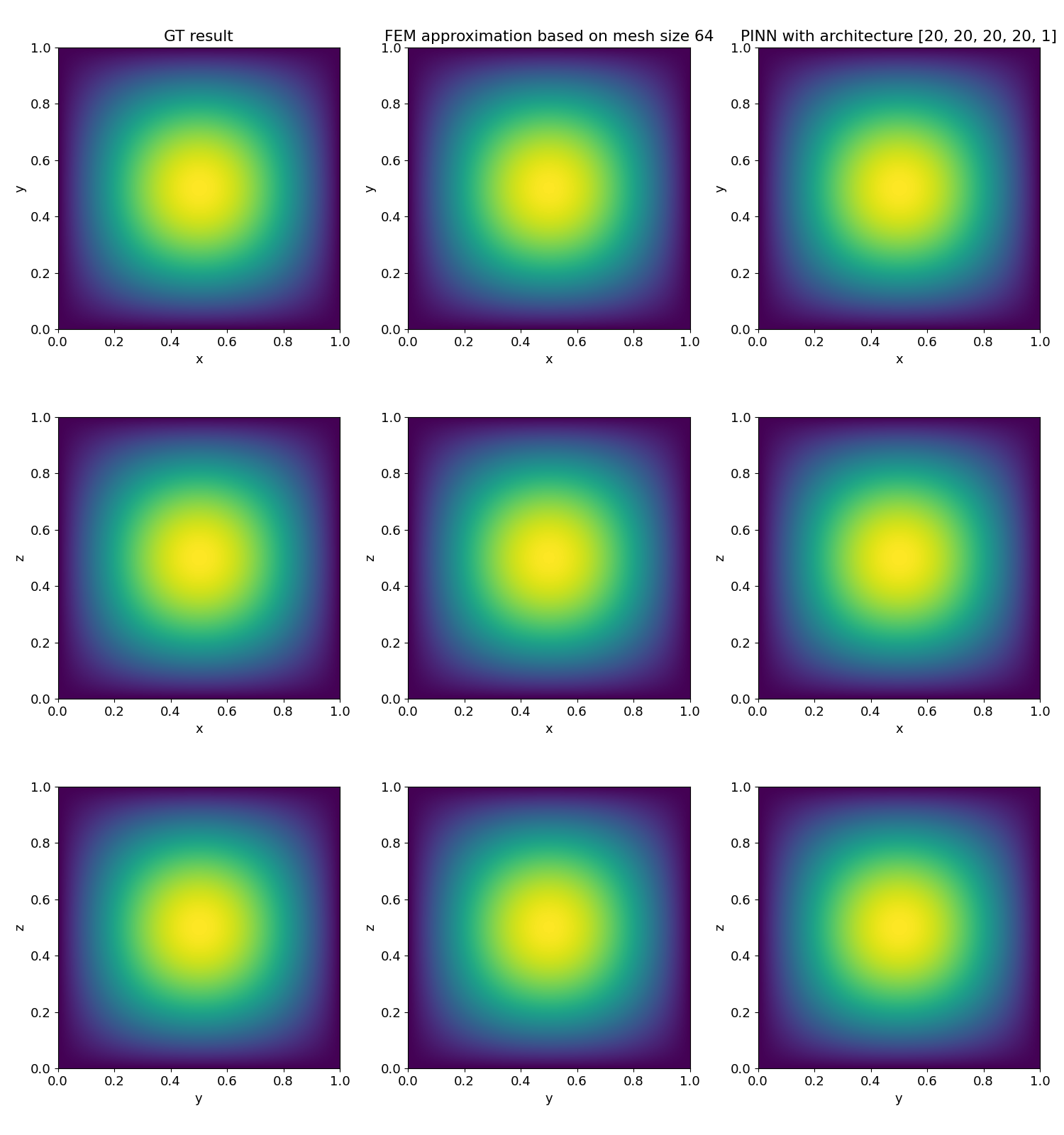}
    \caption{Comparison of the 3D Poisson ground truth solution slices at x,y,z = 0.5 respectively to examples of the FEM and PINN approximations.}
    \label{fig.3DPoisson_solution}
\end{figure*}

\subsubsection*{Results}
\begin{figure*}[ht]
     \centering
     \begin{subfigure}[b]{0.49\textwidth}
         \centering
         \includegraphics[width=\textwidth]{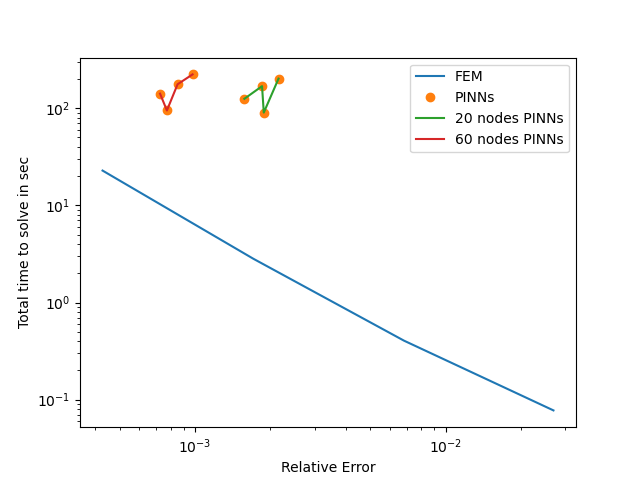}
         \caption{Plot of time to solve FEM and train PINN in sec versus $\ell^2$ relative error.}
         \label{fig.3DPoisson_solve}
     \end{subfigure}
     \hfill
     \begin{subfigure}[b]{0.49\textwidth}
         \centering
         \includegraphics[width=\textwidth]{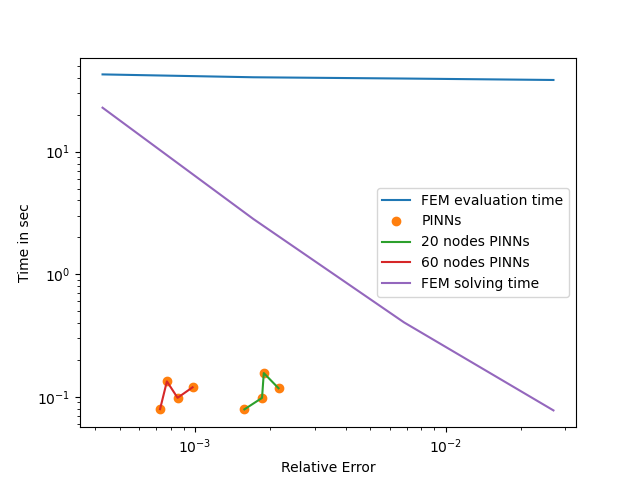}
         \caption{Plot of time to interpolate FEM and evaluate PINN in sec versus relative error. For comparison, the time to solve FEM is also plotted.}
         \label{fig.3DPoisson_test}
     \end{subfigure}
        \caption{Plot for 3D Poisson equation of time in sec versus $\ell^2$ relative error.}
        \label{fig.3DPoisson_timevaccuracy}
\end{figure*}
The PDE is evaluated on a mesh of $150 \times 150 \times 150$ grid points and example results are shown in Figure \ref{fig.3DPoisson_solution}. The time versus accuracy plots can be found in Figure \ref{fig.3DPoisson_timevaccuracy}. As expected, the FEM results with slower computation times have lower relative errors as they are solved on finer meshes as shown in Figure \ref{fig.3DPoisson_solve}. While the PINN approximations are about 1-3 orders of magnitude slower in training time depending on the FEM mesh solution that we compare to, they are able to achieve equal and even higher accuracy in most cases. On the other hand, PINNs outperform FEM when considering the evaluation time as plotted in Figure \ref{fig.3DPoisson_test}. The time to interpolate FEM on a new mesh is 2-3 orders of magnitude slower than the evaluation time of PINNs. Additionally, PINNs are able to achieve equal or higher accuracy scores. If we also take the FEM solving time into account, we find again that this is faster than the evaluation time. However, it is slower than the PINNs evaluation. Therefore, a trained PINN has a lower computation time for the evaluation and a similar to lower relative error as compared to both the FEM solution and evaluation times.

\section{Approximating the Allen-Cahn equation} \label{sec_ACE}
To study the one dimensional Allen-Cahn equation, we consider the following PDE:
\begin{align}\label{eq.1DAllen-Cahn}
    \begin{split}
    \frac{\partial u(t,x)}{\partial t} &=  \epsilon \Delta u - \frac{2}{\epsilon} u(t,x) (1-u(t,x)) (1-2u(t,x)), \quad x \in \Omega = [0,1], \,\, t \in [0,T]\\
    u(t,0) &= u(t,1), \quad t \in [0,T]\\
    u(0,x) &= \frac{1}{2}\left(\frac{1}{2} \sin (x 2 \pi) + \frac{1}{2}\sin(x 16 \pi) \right) + \frac{1}{2}, \quad x \in \Omega,
    \end{split}
\end{align}
where $T = 0.05$ and $\epsilon = 0.01$. As mentioned in Section~\ref{sec_Math_back}, we note that a smaller $\epsilon$ will yields close to piecewise constant solutions, whereas solutions with large $\epsilon$ will be overall more smooth. The $\epsilon$ chosen in this case is due to the inability of PINNs to approximate the Allen-Cahn solution with a smaller $\epsilon$. We have trained the neural network to solve Allen-Cahn with $\epsilon = 0.001$ for various network architectures and with activation functions such as softplus that are typically able to handle discontinuous solutions. However, the PINNs were not able to recover results close to the ground truth solution. We will discuss this case further below. For now, the results refer to $\epsilon = 0.01$.

\subsubsection*{FEM}
As opposed to the Poisson equation, the Allen-Cahn equation has a time dependency that we need to consider in the discretisation. As mentioned in Section \ref{sec.FEM}, we use a semi-implicit Euler strategy to approximate the solution. The weak formulation in a semi-discrete form of the 1D Allen-Cahn equation with Dirichlet $u_\partial = 0$ boundaries can be written as
\begin{align}\label{eq.weak_1DAllen-Cahn}
    \begin{split}
    \int_0^1 &(u_{t+1}(x) - u_{t}(x))v(x) \mathrm{d}x + \epsilon \int_0^1 \langle \nabla u_{t+1}(x), \nabla v(x)\rangle \mathrm{d}x \\
    &+ \frac{2}{\epsilon} \int_0^1 u_{t}(x)(1-u_{t}(x))(1-2u_{t}(x))v(x) \mathrm{d}x = 0,
    \end{split}
\end{align}
for all test functions $v \in H^1_0([0,1])$. The time is discretised with distance of $dt = 1e-3$. For the finite element mesh in space, we choose an interval mesh on the domain $[0,1]$ with varying numbers of cells $n \in \{32,128,512,2048\}$. As in 1D Poisson, we employ piecewise-linear finite element functions ($P_1$). The non-linear variational problem \eqref{eq.weak_1DAllen-Cahn} with periodic boundary conditions is solved using a Newton solver. The FEM approach was implemented using the python toolbox FEniCS.

\subsubsection*{PINNs}
We follow the general PINNs approach and design the loss based on the PDE residual, the boundary residual and the initial condition residual. Additionally, we introduce a weighting of the different terms in the loss functional. We found heuristically this weighting to render the best results.
\begin{align}\label{eq.1DAllen-Cahn_loss}
    \begin{split}
    \text{Loss}(\theta) = \, &\frac{1}{N_f} \sum_{i = 1}^{N_f} \lVert 
    \frac{\partial u_{\theta}(t^i_f,x^i_f)}{\partial t} - \epsilon \Delta u_{\theta} (t^i_f,x^i_f) + \frac{2}{\epsilon} u_{\theta}(t^i_f,x^i_f) (1-u_{\theta}(t^i_f,x^i_f)) (1-2u_{\theta}(t^i_f,x^i_f)) \rVert^2_2 \\
    &+ \frac{1}{N_g} \sum_{k = 1}^{N_g}\lVert u_{\theta}(t^j_g,0) - u_{\theta}(t^j_g,1) \rVert^2_2 \\
    &+ \frac{1000}{N_h} \sum_{k = 1}^{N_h}\lVert u_{\theta}(0,x^k_h) -  \frac{1}{2}\left(\frac{1}{2} \sin (2 \pi x^k_h) + \frac{1}{2}\sin(16 \pi x^k_h) \right) + \frac{1}{2}\rVert^2_2,
    \end{split}
\end{align}
with $u_{\theta}$ the neural network and $\theta$ the trained weights. We train the network on $N_f = 20,000$ collocation points $(t_f^i,x_f^i) \in [0,0.05]\times [0,1]$ that are sampled using Latin Hybercubes. The training points for the boundary with $N_g = 250$ and for the initial condition with $N_h = 500$ are also sampled in each epoch with Latin Hybercube sampling. The network architecture is a feed-forward dense neural network with a $\tanh$ activation functional. We consider architectures of the following 14 different sizes: [20,20,20,1], [100,100,100,1], [500,500,500,1], [20,20,20,20,1], [100,100,100,100,1], [500,500,500,500,1], [20,20,20,20,20,1], [100,100,100,100,100,1],  [500,500,500,500,500,1], [20,20,20,20,20,20,1], [100,100,100,100,100,100,1], [500,500,500,500,500,500,1], [20,20,20,20,20,20,20,1] and [100,100,100,100,100,100,100,1]. For a network of size [500,500,500,500,500,500,500,1] we run out of memory on the GPU available to us. This constitutes the limitation of our training. For the optimisation we first run the Adam optimiser with learning rate $1e-4$ for 7000 epochs over the initial loss alone. We have found the network struggling to learn the initial condition when the optimisation is run on the full loss directly. After the 7000 epochs optimising the initial loss, we run the Adam optimiser for $50,000$ epochs on the full loss function \eqref{eq.1DAllen-Cahn_loss}. Finally, we refine the optimisation using L-BFGS. 
\begin{figure*}[t]
    \centering
    \includegraphics[width=\textwidth]{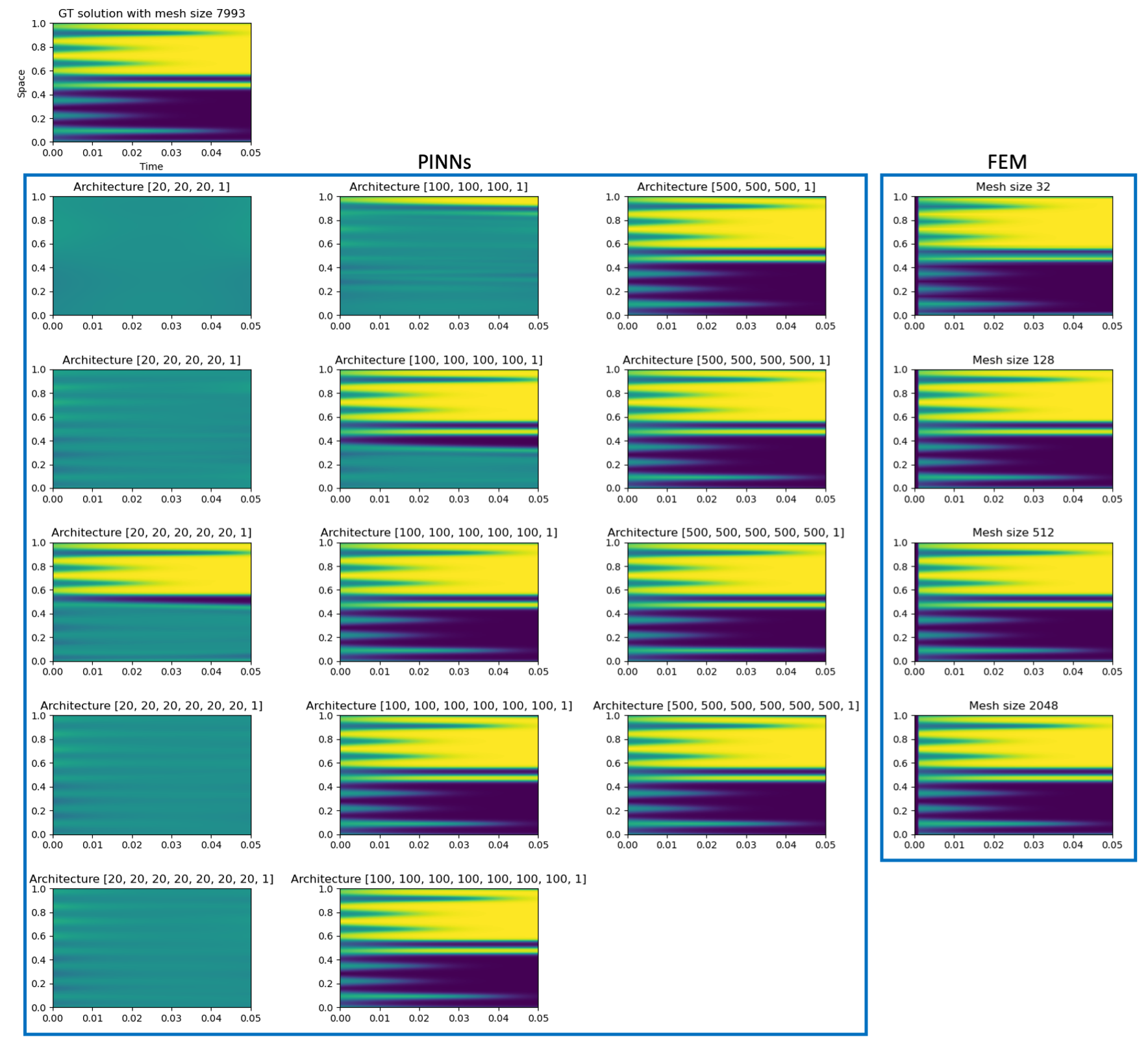}
    \caption{Comparison of the 1D Allen-Cahn solution approximated by FEM and PINNs.}
    \label{fig.1DAllen-Cahn_solution}
\end{figure*}

\subsubsection*{Results}
The PDE approximations with FEM and PINNs are compared on the mesh  of the ground truth solution spanning $[0,1]$ with 7993 mesh points and at a time discretisation of $\frac{1}{3} e-4$. The fine-meshed FEM approximation for the ground truth solution was derived using implicit Euler. The FEM and the PINN approximations compared to the ground truth solution are displayed in Figure \ref{fig.1DAllen-Cahn_solution} for the different mesh sizes and network architectures. FEM is able to recover the solution of the PDE at all mesh sizes. However, closer inspection of the result for a mesh of size $32$ shows slight errors along the diffusive interface. On the other hand, the ability of a PINN to approximate the PDE solution well is dependent on the architecture and the number of free parameters or weights that are to be determined. While all architectures with 20 nodes (cf. Figure \ref{fig.1DAllen-Cahn_solution} column 1) are not able to recover the solution whatsoever, networks with 100 nodes per layer are able to be trained for the solution. We can clearly observe a progression based on the number of layers. Finally, neural networks that have 500 nodes per layer are all able to approximate the solution well. 

\begin{figure*}[ht]
     \centering
     \begin{subfigure}[b]{0.49\textwidth}
         \centering
         \includegraphics[width=\textwidth]{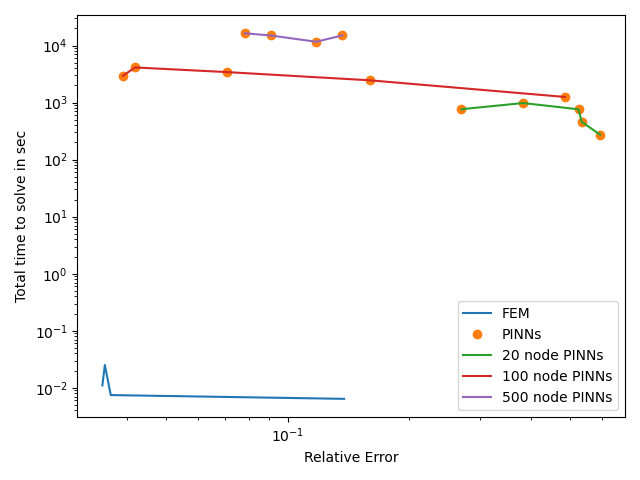}
         \caption{Plot of time to solve FEM and train PINN in sec versus $\ell^2$ relative error.}
         \label{fig.Allen-Cahn_solve}
     \end{subfigure}
     \hfill
     \begin{subfigure}[b]{0.49\textwidth}
         \centering
         \includegraphics[width=\textwidth]{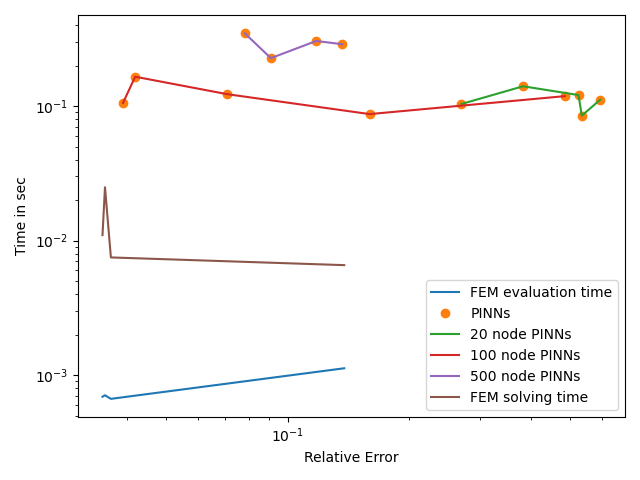}
         \caption{Plot of time to interpolate FEM and evaluate PINN in sec versus $\ell^2$ relative error. For comparison, the time to solve FEM is plotted.}
         \label{fig.Allen-Cahn_test}
     \end{subfigure}
        \caption{Plot for 1D Allen-Cahn equation of time in sec versus $\ell^2$ relative error.}
        \label{fig.1DAllen-Cahn_plots}
\end{figure*}
Let us compare the solution and evaluation time versus the accuracy for the FEM and PINN approximations as shown in Figure \ref{fig.1DAllen-Cahn_plots}. For the solution time, that is the time to solve the PDE using FEM and the time to train the neural network in PINN, FEM  is 5-6 orders of magnitude faster than PINNs. This is mainly due to the size of the neural networks. Here, the complexity of the PDE requires larger architectures to be able to capture the solution. Some of the PINN architectures are then able though to achieve similar relative errors to FEM. This is especially true for networks with 100 or 500 nodes per layer as also seen in Figure \ref{fig.1DAllen-Cahn_solution}. The evaluation time of the PINN is plotted against the evaluation time of FEM and the solution time of FEM in Figure \ref{fig.Allen-Cahn_test}. The calculation time advantage of FEM evaluation drops to about one order of magnitude as compared to the solution time. While it is to be expected that the evalutation of a neural network is much faster then the training, FEM is still faster in both solution and evaluation time.

\begin{figure*}[t]
    \centering
    \includegraphics[width=0.3\textwidth]{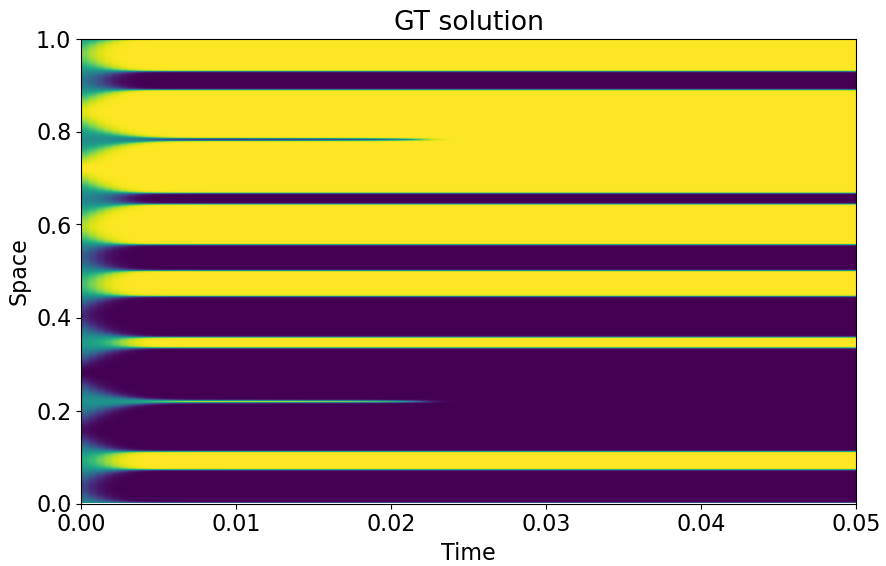}
    \caption{FEM solution of the 1D Allen-Cahn equation for $\epsilon = 0.001$ derived on a mesh with 7993 cells.}
    \label{fig.1DAllen-Cahn_small_eps}
\end{figure*}

Compared to the other PDEs that we are considering in this study, the Allen-Cahn equation needs slight modification, i.e. weighting of the loss and pre-training on only part of the loss functional. This is due to the difficulty we have found the network to have in learning the PDE solution. In addition, we should note that we have also attempted to train a PINN for Allen-Cahn \eqref{eq.1DAllen-Cahn} with $\epsilon=0.001$. The resulting solution  -- shown in Figure \ref{fig.1DAllen-Cahn_small_eps} -- becomes close to binary after a certain amount of time. We had discussed this effect in Section~\ref{sec_Math_back}. This renders very large gradients or  discontinuities in the solution. We trained PINNs with different activation functions such as softplus or ReLU that are typically able to handle discontinuities. However, all results were insufficient to be considered an approximation. This speaks to the assumption that PINNs in the vanilla form are not well equipped to handle discontinuous solutions; this may also be due to PINNs solving the strong PDE, rather than a weak form. Variations of the vanilla PINNs approach might be able to obtain satisfactory approximations. However, as we are only considering the vanilla approach, this goes beyond the scope of the paper. We make a note that FEM is able to approximate Allen-Cahn with $\epsilon=0.001$ albeit we should anticipate the use of a finer mesh to accurate represent the diffuse interface. 

\section{Approximating the Semilinear Schrödinger equation} \label{Sec_SSE}
Finally, we investigate the semilinear Schrödinger equation in one and two space-dimensions. The specifications of the one-dimensional case are taken from the original PINNs paper of Raissi et al. \cite{Raissi2019}. Note that the semilinear Schrödinger equation has a complex-valued solution; further increasing the difficulty of its approximation.
\subsection{1D}
\begin{figure*}[t]
    \centering
    \includegraphics[width=\textwidth]{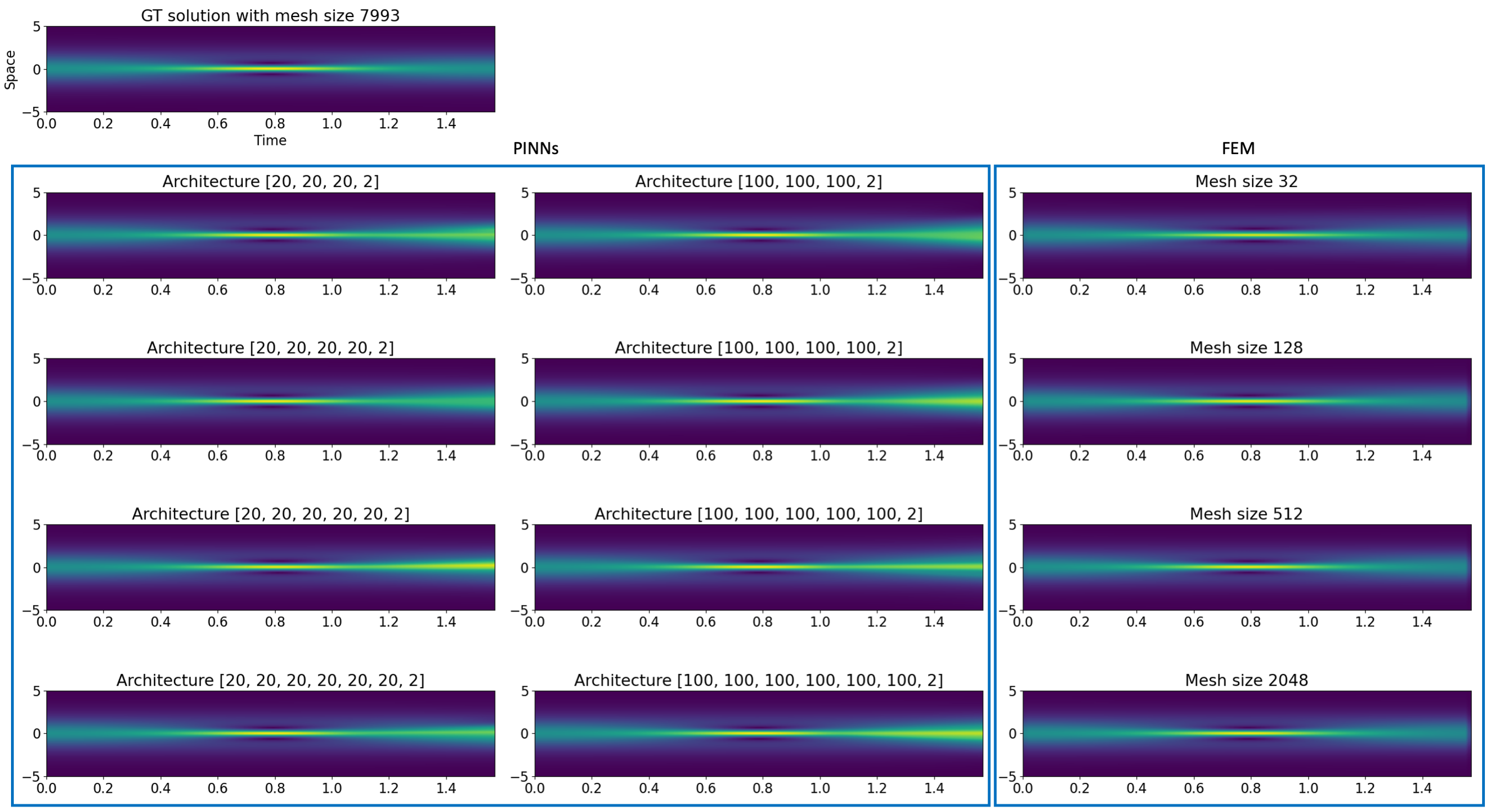}
    \caption{Comparison of the 1D semilinear Schrödinger solution $\lvert h(t,x) \rvert = \sqrt{u_R^2(t,x) + u_I^2(t,x)}$ approximated by FEM and PINNs of different mesh and architecture sizes.}
    \label{fig.1DSchroedinger_solution}
\end{figure*}

In the one-dimensional case, we consider the semilinear Schrödinger equation with periodic boundary conditions, as follows:
\begin{align*}
    \mathrm{i}\frac{\partial h(t,x)}{\partial t} &= -0.5 \Delta h(t,x) - \lvert h(t,x) \rvert^2 h(t,x) &x \in [-5,5], \,\, t \in [0,\pi /2],\\
    h(0,x) &= 2 \text{sech}(x),  &x \in [-5,5],\\
    h(t,-5) &= h(t,5), &t \in [0,\pi /2],\\
    \frac{\partial h(t,-5)}{\partial x} &= \frac{\partial h(t,5)}{\partial x},  &t \in [0,\pi /2].
\end{align*}
We note that the identical problem had also been considered by \cite{Raissi2019}. As $h(t,x)$ is a complex valued function, the semilinear Schrödinger equation is solved for $h(t,x) =: u_R(t,x) + \mathrm{i}u_I(t,x)$, with $u_R(t,x)$ the real part and $u_I(t,x)$ the imaginary part. Example results are visualised for $\lvert h(t,x) \rvert = \sqrt{u_R^2(t,x) + u_I^2(t,x)}$ in Figure \ref{fig.1DSchroedinger_solution}. 

\subsubsection*{FEM}
Similar to the 1D Allen-Cahn equation that we considered, we use a semi-implicit Euler strategy to approximate the 1D Schrödinger equation in time. We can then write the weak formulation for the real and imaginary parts of the PDE, assuming Dirichlet $u_\partial = 0$ boundaries, separately as
\begin{align*}
    \int_{-5}^5 (u^I_{t+1}(x) - u^I_{t}(x))v^R(x) \mathrm{d}x - 0.5 \int_{-5}^5 \langle \nabla u^R_{t+1}(x), \nabla v^R(x)\rangle \mathrm{d}x - \lvert h_{t}(x) \rvert^2 \int_{-5}^5 u^R_{t+1}(x)v^R(x) \mathrm{d}x &= 0,\\
    \int_{-5}^5 (u^R_{t+1}(x) - u^R_{t}(x))v^I(x) \mathrm{d}x + 0.5 \int_{-5}^5 \langle \nabla u^I_{t+1}(x), \nabla v^I(x)\rangle \mathrm{d}x + \lvert h_{t}(x) \rvert^2 \int_{-5}^5 u^I_{t+1}(x)v^I(x) \mathrm{d}x &= 0,
\end{align*}
for all test functions $v^R, v^I \in H_0^1(\Omega)$. The time is discretised with $dt = 1e-3$ and we define the finite elements on an interval mesh in $[-5,5]$. We consider meshes with $n \in \{32,128,512,2048\}$ cells. We again employ piecewise linear finite element basis functions on those cells ($P_1$) and use the generalised minimal residual method (gmres) for a solver.

\subsubsection*{PINNs}
Let us now define the loss functional used in the neural network approximation of the PDE solution. Again, we consider the residual of the PDE, inital condition and boundary conditions as follows:
\begin{align*}\label{eq.1DSchroedinger_loss}
    \begin{split}
    \text{Loss}(\theta) = \, &\frac{1}{N_f} \sum_{i = 1}^{N_f} (\lVert \frac{\partial u^I_{\theta}(t^i_f,x^i_f)}{\partial t} - \epsilon \Delta u^R_{\theta} (t^i_f,x^i_f) -  \lvert h_{\theta}(t^i_f,x^i_f) \rvert^2  u^R_{\theta}(t^i_f,x^i_f) \rVert^2_2 \\
    &\qquad + \lVert \frac{\partial u^R_{\theta}(t^i_f,x^i_f)}{\partial t} + \epsilon \Delta u^I_{\theta} (t^i_f,x^i_f) + \lvert h_{\theta}(t^i_f,x^i_f) \rvert^2 u^I_{\theta}(t^i_f,x^i_f) \rVert^2_2 )\\
    &+ \frac{1}{N_g} \sum_{k = 1}^{N_g}\left(\lVert u^R_{\theta}(t^j_g,-5) - u^R_{\theta}(t^j_g,5) \rVert^2_2 + \lVert u^I_{\theta}(t^j_g,-5) - u^I_{\theta}(t^j_g,5) \rVert^2_2 \right)\\
    &+ \frac{1}{N_g} \sum_{k = 1}^{N_g}\left(\lVert \frac{\partial u^R_{\theta}(t^j_g,-5)}{\partial x} - \frac{\partial u^R_{\theta}(t^j_g,5)}{\partial x} \rVert^2_2 + \lVert \frac{\partial u^I_{\theta}(t^j_g,-5)}{\partial x} - \frac{\partial u^I_{\theta}(t^j_g,5)}{\partial x} \rVert^2_2 \right)\\
    &+ \frac{1}{N_h} \sum_{k = 1}^{N_h}\left(\lVert u^R_{\theta}(0,x^k_h) -  2 \text{sech}(x^k_h)\rVert^2_2 + \lVert u^I_{\theta}(0,x^k_h)\rVert^2_2\right),
    \end{split}
\end{align*}
where $h_{\theta}(t,x) = [u^R_{\theta}(t,x),u^I_{\theta}(t,x)]$ the neural network the produces the real and imaginary parts of the PDE solution of the network output with weights $\theta$. The network is trained on $N_f = 20,000$ collocation points $(t^i_f,x^i_f) \in [0,\pi /2] \times [-5,5]$ and with $N_g = 50$ point on the boundary and $N_h = 50$ points for the initial condition using Latin Hybercube sampling. The network architecture is a feed-forward dense neural network with $\tanh$ activation function. In these specifications, we have followed the original vanilla PINNs paper \cite{Raissi2019}. We investigate the performance of 8 network architectures, that are: [20,20,20,2], [100,100,100,2], [20,20,20,20,2], [100,100,100,100,2], [20,20,20,20,20,2], [100,100,100,100,100,2], [20,20,20,20,20,20,2], and [100,100,100,100,100,100,2]. We employ the Adam optimiser for $50,000$ epochs and a learning rate of $1e-4$. Afterwards, we use the L-BFGS optimiser to refine the training results.

\subsubsection*{Results}
\begin{figure*}[t]
     \centering
     \begin{subfigure}[b]{0.32\textwidth}
         \centering
         \includegraphics[width=\textwidth]{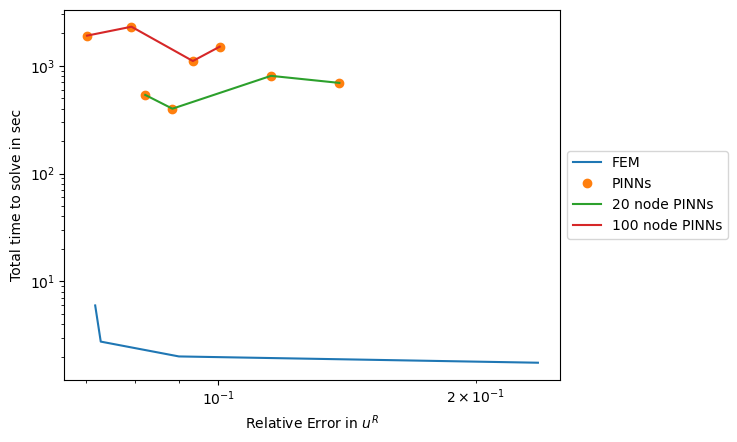}
         \caption{Solution time versus accuracy for the real part of the PDE $u^R$.}
         \label{fig.1DSchroedinger_solve_r}
     \end{subfigure}
     \hfill
     \begin{subfigure}[b]{0.32\textwidth}
         \centering
         \includegraphics[width=\textwidth]{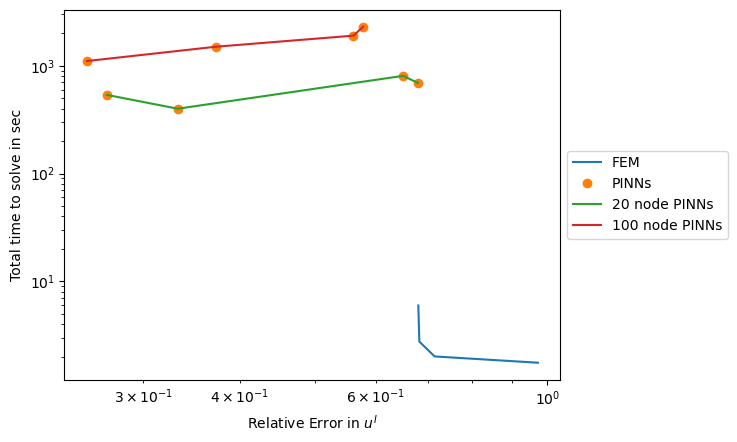}
         \caption{Solution time versus accuracy for the imaginary part of the PDE $u^I$.}
         \label{fig.1DSchroedinger_solve_i}
     \end{subfigure}
     \hfill
     \begin{subfigure}[b]{0.32\textwidth}
         \centering
         \includegraphics[width=\textwidth]{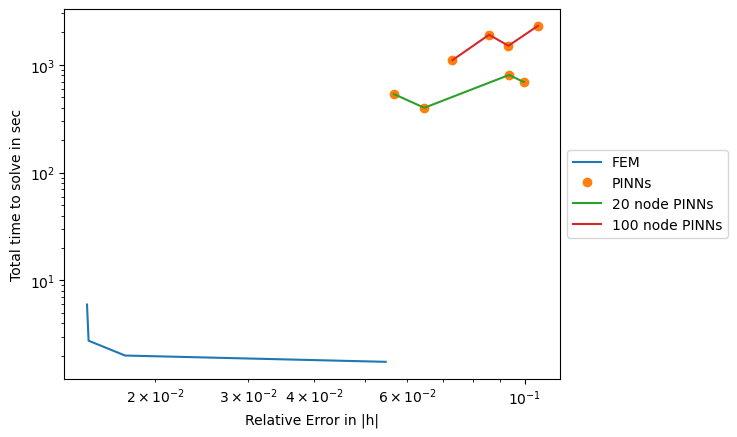}
         \caption{Solution time versus accuracy for the modulus $\lvert h \rvert$.}
         \label{fig.1DSchroedinger_solve_h}
     \end{subfigure}
     \hfill
     \begin{subfigure}[b]{0.32\textwidth}
         \centering
         \includegraphics[width=\textwidth]{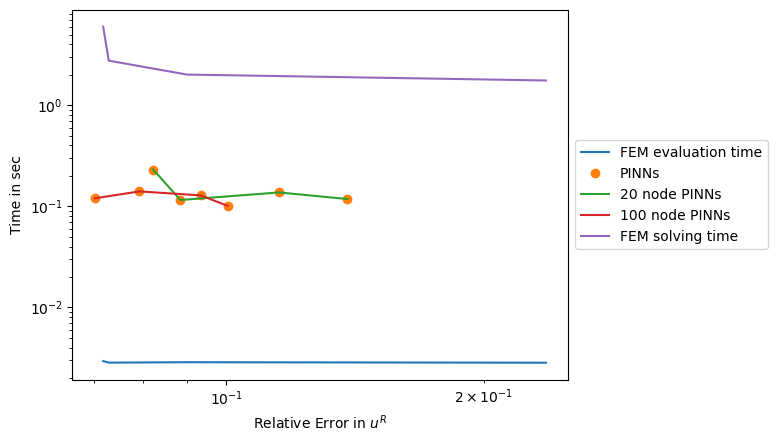}
         \caption{Evaluation time versus accuracy for the real part of the PDE $u^R$. For comparison, the time to solve FEM is plotted.}
         \label{fig.1DSchroedinger_test_r}
     \end{subfigure}
     \hfill
     \begin{subfigure}[b]{0.32\textwidth}
         \centering
         \includegraphics[width=\textwidth]{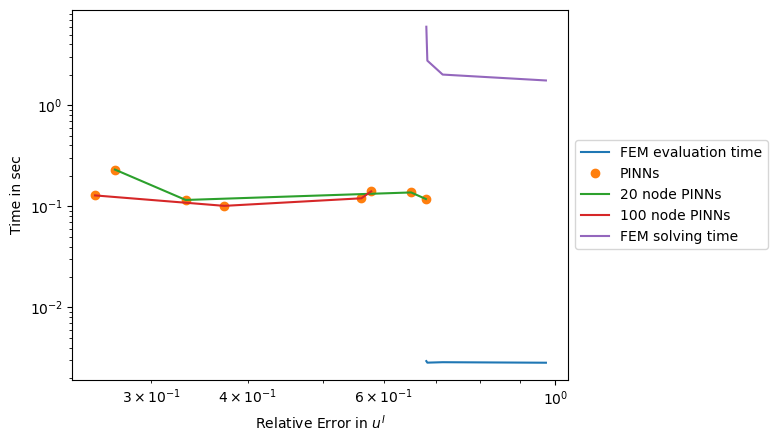}
         \caption{Evaluation time versus accuracy for the imaginary part of the PDE $u^I$. For comparison, the time to solve FEM is plotted.}
         \label{fig.1DSchroedinger_test_i}
     \end{subfigure}
     \hfill
     \begin{subfigure}[b]{0.32\textwidth}
         \centering
         \includegraphics[width=\textwidth]{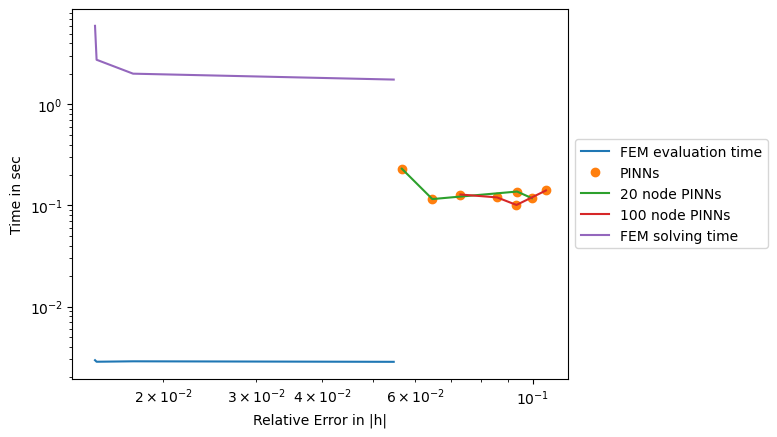}
         \caption{Evaluation time versus accuracy for the modulus $\lvert h \rvert$. For comparison, the time to solve FEM is plotted.}
         \label{fig.1DSchroedinger_test_h}
     \end{subfigure}
        \caption{Plot for 1D Schrödinger equation of time in sec versus $\ell^2$ relative error. Plots are split for the real and imaginary parts of the PDE as well as $\lvert h \rvert = \sqrt{u^2_R + u_I^2}$.}
        \label{fig.1DSchroedinger_plots}
\end{figure*}
The results of the FEM and PINNs approximation are compared on the mesh of the ground truth solution that has a size of 7993 cells and a time discretisation with $dt = \frac{1}{3}e-4$. The resulting approximations for $\lvert h(t,x) \rvert = \sqrt{u^2_R(t,x) + u_I^2(t,x)}$ are displayed in Figure \ref{fig.1DSchroedinger_solution}. Already visually, we notice that the PINN approximations become slightly less accurate for larger time instance than the FEM approximations. Quantitatively speaking, the time versus accuracy plots give a broader overview of the performance of each method. The plots are shown in Figure \ref{fig.1DSchroedinger_plots}. Focusing on the modulus $\lvert h \rvert$, FEM both has a lower solution time by 2 orders of magnitude and a lower relative error than any neural network approximation as shown in Figure \ref{fig.1DSchroedinger_solve_h}. Considering the evaluation time for both methods alone, FEM continues to outperform PINNs in time and accuracy. However, if we compare the FEM solving time to the PINNs evaluation time, neural networks are able to evaluate solutions faster by one order of magnitude. Nonetheless, FEM remains to produce results with higher accuracy for $|h|$. The results are inconclusive for $u_R$ and $u_I$: the FEM solution is considerably faster but less or hardly more accurate than the PINNs solution as shown in Figures \ref{fig.1DSchroedinger_test_r} and \ref{fig.1DSchroedinger_test_i}. 

\subsection{2D}
\begin{figure*}[t]
    \centering
    \includegraphics[width=\textwidth]{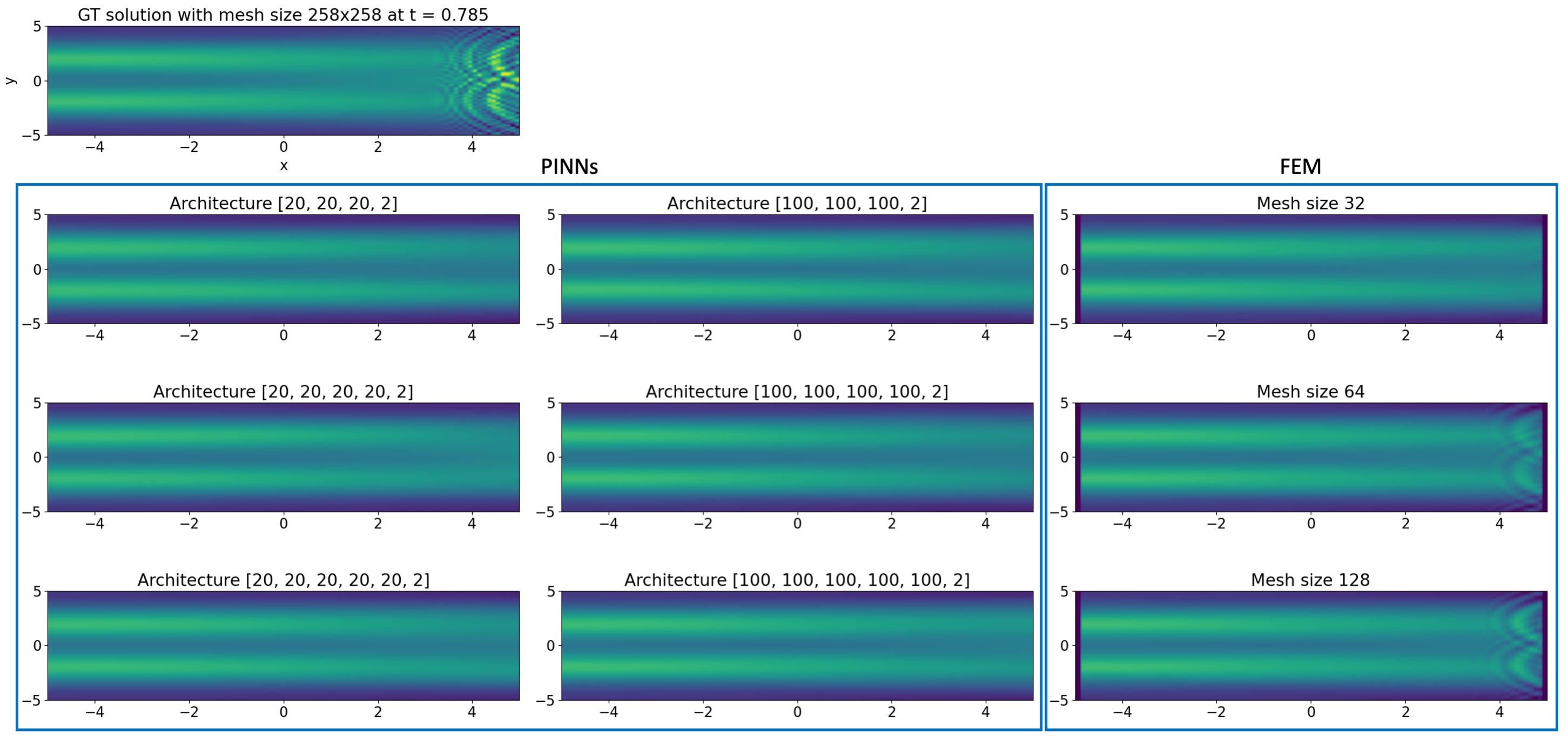}
    \caption{Comparison of the 2D semilinear Schrödinger solution $\lvert h(t,x) \rvert = \sqrt{u_R^2(t,x) + u_I^2(t,x)}$ at time $t = \pi/4$ approximated by FEM and PINNs of different mesh and architecture sizes.}
    \label{fig.2DSchroedinger_solution}
\end{figure*}

Let us finally move to the two-dimensional semilinear Schrödinger equation. We consider periodic boundary conditions and define the initial condition as follows:
\begin{align*}
    \mathrm{i}\frac{\partial h(t,x,y)}{\partial t} &= -0.5 \Delta h(t,x,y) - \lvert h(t,x,y) \rvert^2 h(t,x,y), &x,y \in [-5,5], \,\, t \in [0,\pi /2]\\
    h(0,x,y) &= \text{sech}(x)+0.5\text{sech}(y-2)+0.5\text{sech}(y+2), &x,y \in [-5,5]\\
    h(t,-5,y) &= h(t,5,y), &t \in [0,\pi /2],\,\, y\in [-5,5]\\
    h(t,x,-5) &= h(t,x,5), &t \in [0,\pi /2],\,\, x\in [-5,5]\\
    \frac{\partial h(t,-5,y)}{\partial x} &= \frac{\partial h(t,5,y)}{\partial x}, &t \in [0,\pi /2],\,\, y\in [-5,5]\\
    \frac{\partial h(t,x,-5)}{\partial y} &= \frac{\partial h(t,x,5)}{\partial y} &t \in [0,\pi /2],\,\, x\in [-5,5].
\end{align*}
The PDE is complex-valued and we approximate the solution for $h(t,x,y) =: u_R(t,x,y) + \mathrm{i}u_I(t,x,y)$. The approximate solutions for the modulus for $\lvert h(t,x,y) \rvert = \sqrt{u_R^2(t,x,y) + u_I^2(t,x,y)}$ are shown in Figure \ref{fig.2DSchroedinger_solution}.

\subsubsection*{FEM}
Similar to the one-dimensional case of the Schrödinger equation, we split the weak formulation of the PDE (again for Dirichlet $u_\partial = 0$ boundaries)  into its real and imaginary parts:
\begin{align*}
    \int_{-5}^5\int_{-5}^5 (u^I_{t+1}(x,y) &- u^I_{t}(x,y))v^R(x,y) \mathrm{d}x \mathrm{d}y- 0.5 \int_{-5}^5\int_{-5}^5 \langle \nabla u^R_{t+1}(x,y), \nabla v^R(x,y)\rangle \mathrm{d}x \mathrm{d}y \\
    &- \lvert h_{t}(x,y) \rvert^2 \int_{-5}^5\int_{-5}^5 u^R_{t+1}(x,y)v^R(x,y) \mathrm{d}x \mathrm{d}y = 0,\\
    \int_{-5}^5\int_{-5}^5 (u^R_{t+1}(x,y) &- u^R_{t}(x,y))v^I(x,y) \mathrm{d}x \mathrm{d}y + 0.5 \int_{-5}^5\int_{-5}^5 \langle \nabla u^I_{t+1}(x,y), \nabla v^I(x,y)\rangle \mathrm{d}x \mathrm{d}y \\
    &+ \lvert h_{t}(x,y) \rvert^2 \int_{-5}^5\int_{-5}^5 u^I_{t+1}(x,y)v^I(x) \mathrm{d}x \mathrm{d}y = 0,
\end{align*}
for all test functions $v^R,v^I \in H_0^{\Omega}$. Time is discretised with $dt = 1e-3$ and the finite elements are defined on a rectangular mesh in the domain $[-5,5] \times [-5,5]$. We consider meshes with $\{(16,16),(32,32),(40,40),(64,64),(128,128)\}$ squares, each again being split into two triangles on which we define piecewise-linear finite element basis functions ($P_1$). We again use gmres to solve the linear system.

\subsubsection*{PINNs}
The loss functional for the neural network approximations are defined for both the real and imaginary parts of the PDE:
\begin{align*}
    \text{Loss}(\theta) = \, &\frac{1}{N_f} \sum_{i = 1}^{N_f} (\lVert \frac{\partial u^I_{\theta}(t^i_f,x^i_f,y^i_f)}{\partial t} - \epsilon \Delta u^R_{\theta} (t^i_f,x^i_f,y^i_f) -  \lvert h_{\theta}(t^i_f,x^i_f,y^i_f) \rvert^2  u^R_{\theta}(t^i_f,x^i_f,y^i_f) \rVert^2_2 \\
    &\qquad + \lVert \frac{\partial u^R_{\theta}(t^i_f,x^i_f,y^i_f)}{\partial t} + \epsilon \Delta u^I_{\theta} (t^i_f,x^i_f,y^i_f) + \lvert h_{\theta}(t^i_f,x^i_f,y^i_f) \rvert^2 u^I_{\theta}(t^i_f,x^i_f,y^i_f) \rVert^2_2 )\\
    &+ \frac{1}{N_g} \sum_{k = 1}^{N_g}(\lVert u^R_{\theta}(t^j_g,-5,y^j_g) - u^R_{\theta}(t^j_g,5,y^j_g) \rVert^2_2 + \lVert u^I_{\theta}(t^j_g,-5,y^j_g) - u^I_{\theta}(t^j_g,5,y^j_g) \rVert^2_2\\
    &\qquad + \lVert u^R_{\theta}(t^j_g,x^j_g,-5) - u^R_{\theta}(t^j_g,x^j_g,5) \rVert^2_2 + \lVert u^I_{\theta}(t^j_g,x^j_g,-5) - u^I_{\theta}(t^j_g,x^j_g,5) \rVert^2_2 )\\
    &+ \frac{1}{N_g} \sum_{k = 1}^{N_g}\left(\lVert \frac{\partial u^R_{\theta}(t^j_g,-5,y^j_g)}{\partial x} - \frac{\partial u^R_{\theta}(t^j_g,5,y^j_g)}{\partial x} \rVert^2_2 + \lVert \frac{\partial u^I_{\theta}(t^j_g,-5,y^j_g)}{\partial x} - \frac{\partial u^I_{\theta}(t^j_g,5,y^j_g)}{\partial x} \rVert^2_2 \right)\\
    &+ \frac{1}{N_g} \sum_{k = 1}^{N_g}\left(\lVert \frac{\partial u^R_{\theta}(t^j_g,x^j_g,-5)}{\partial y} - \frac{\partial u^R_{\theta}(t^j_g,x^j_g,5)}{\partial y} \rVert^2_2 + \lVert \frac{\partial u^I_{\theta}(t^j_g,x^j_g,-5)}{\partial y} - \frac{\partial u^I_{\theta}(t^j_g,x^j_g,5)}{\partial y} \rVert^2_2 \right)\\
    &+ \frac{1}{N_h} \sum_{k = 1}^{N_h}\left(\lVert u^R_{\theta}(0,x^k_h,y^k_h) -  \text{sech}(x^k_h) - 0.5\text{sech}(y^k_h-2) - 0.5\text{sech}(y^k_h+2) \rVert^2_2\right),
\end{align*}
for $h_{\theta}(t,x,y) = u^R_{\theta}(t,x,y) + \mathrm{i}u^I_{\theta}(t,x,y)$ the neural network with an output size 2 for the real and imaginary parts and $\theta$ the network weights that are determined by training. The loss functional is optimised based on latin hybercube sampled collocation points $(t_f^i,x_f^i,y_f^i) \in [0, \pi /2] \times [-5,5] \times [-5,5]$ with $N_f = 5,000$ points for the domain, $N_g = 100$ points for the boundary and $N_h = 100$ points for the inital condition. We chose feed-forward dense neural netwrok in the architecture design with $\tanh$ the activation function. The results are compared for 6 different neural network sizes: [20,20,20,2], [100,100,100,2], [20,20,20,20,2], [100,100,100,100,2], [20,20,20,20,20,2], and [100,100,100,100,100,2]. The Adam optimiser is run for $50,000$ epochs with a learning rate of $1e-3$ before employing the L-BFGS optimiser.

\begin{figure*}[t]
     \centering
     \begin{subfigure}[b]{0.32\textwidth}
         \centering
         \includegraphics[width=\textwidth]{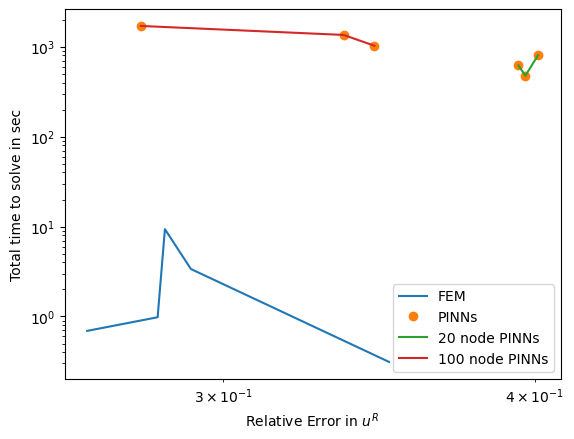}
         \caption{Solution time versus accuracy for the real part of the PDE $u^R$.}
         \label{fig.2DSchroedinger_solve_r}
     \end{subfigure}
     \hfill
     \begin{subfigure}[b]{0.32\textwidth}
         \centering
         \includegraphics[width=\textwidth]{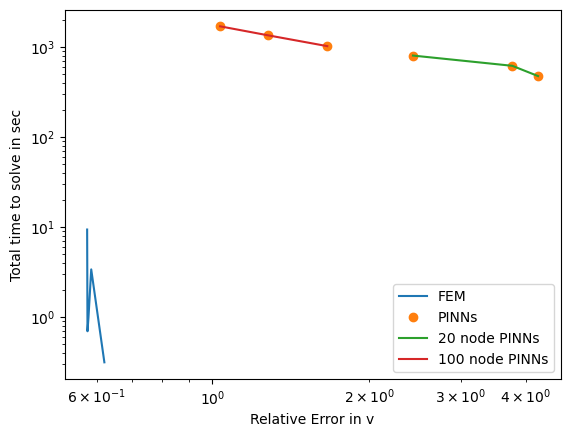}
         \caption{Solution time versus accuracy for the imaginary part of the PDE $u^I$.}
         \label{fig.2DSchroedinger_solve_i}
     \end{subfigure}
     \hfill
     \begin{subfigure}[b]{0.32\textwidth}
         \centering
         \includegraphics[width=\textwidth]{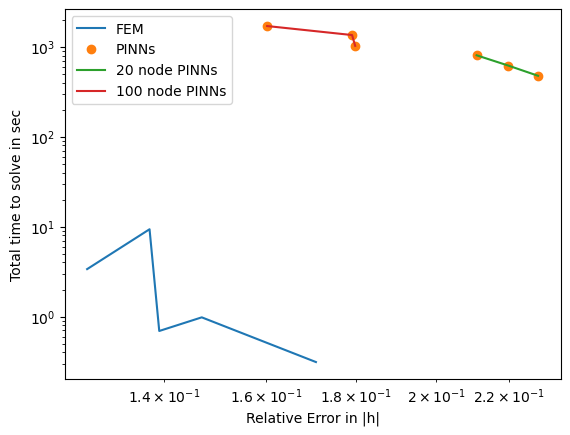}
         \caption{Solution time versus accuracy for the modulus $\lvert h \rvert$.}
         \label{fig.2DSchroedinger_solve_h}
     \end{subfigure}
     \hfill
     \begin{subfigure}[b]{0.32\textwidth}
         \centering
         \includegraphics[width=\textwidth]{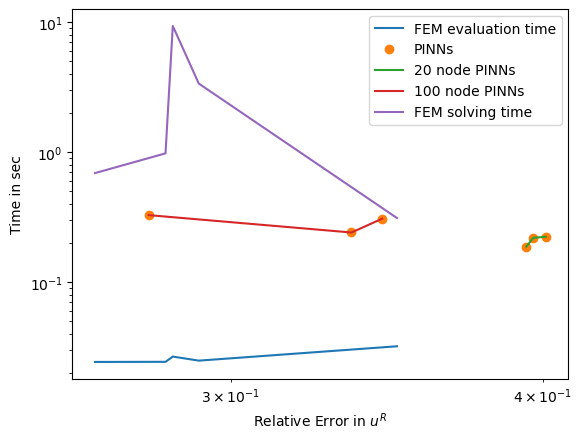}
         \caption{Evaluation time versus accuracy for the real part of the PDE $u^R$. For comparison, the time to solve FEM is plotted.}
         \label{fig.2DSchroedinger_test_r}
     \end{subfigure}
     \hfill
     \begin{subfigure}[b]{0.32\textwidth}
         \centering
         \includegraphics[width=\textwidth]{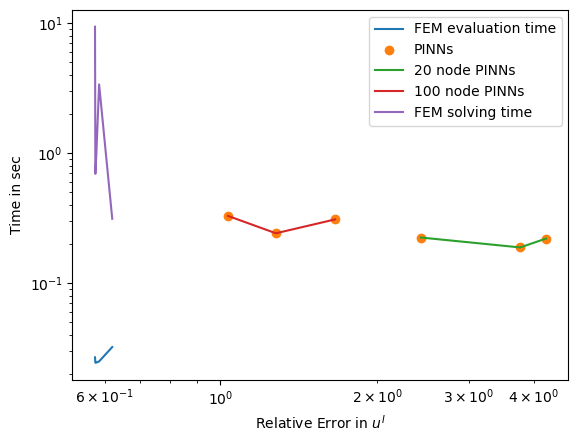}
         \caption{Evaluation time versus accuracy for the imaginary part of the PDE $u^I$. For comparison, the time to solve FEM is plotted.}
         \label{fig.2DSchroedinger_test_i}
     \end{subfigure}
     \hfill
     \begin{subfigure}[b]{0.32\textwidth}
         \centering
         \includegraphics[width=\textwidth]{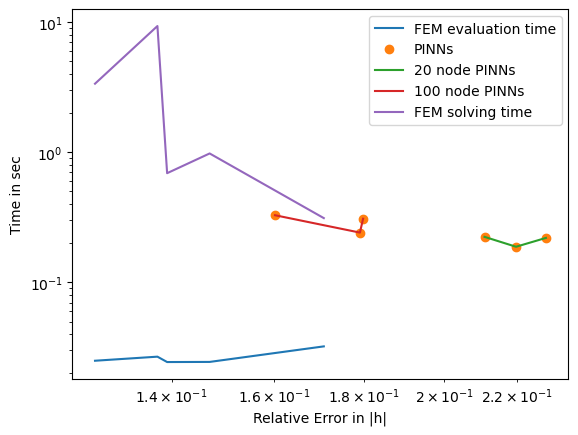}
         \caption{Evaluation time versus accuracy for modulus $\vert h \rvert$. For comparison, the time to solve FEM is plotted.}
         \label{fig.2DSchroedinger_test_h}
     \end{subfigure}
        \caption{Plot for 2D Schrödinger equation of time in sec versus $\ell^2$ relative error. Plots are split for the real and imaginary parts of the PDE as well as $\lvert h \rvert = \sqrt{u^2_R + u_I^2}$.}
        \label{fig.2DSchroedinger_plots}
\end{figure*}

\subsubsection*{Results}
A first look at the visualisation of the results at an example time $t = \pi/4$ in Figure \ref{fig.2DSchroedinger_solution} shows that PINNs is having difficulties to recover the wave-type shapes at the right hand side boundary of $x$. While PINNs are able to reproduce the main features of the PDE solution, the examples shown lack the detailed features containing edges and discontinuities. For FEM, we observe that coarse meshes similarly fail at recovering the fine details. However, with finer meshes, FEM is able to correctly solve the PDE. Considering the quantified results on time versus accuracy in Figure \ref{fig.2DSchroedinger_plots}, the difference in accuracy between PINNs and FEM become less prudent in the real part. However, for both solving time in Figures \ref{fig.2DSchroedinger_solve_r}-\ref{fig.1DSchroedinger_solve_h} and evaluation time in Figures \ref{fig.2DSchroedinger_test_r}-\ref{fig.1DSchroedinger_test_h} FEM clearly outperforms PINNs by 2-3 and 1 order of magnitude, respectively. Taking the FEM solution time into account, we can see that PINNs slightly outperforms in their evaluation time with similar relative error. Let us also note that the PINN approximation of the imaginary part of the complex-valued PDE solution is significantly less accurate compared to the FE method.

\section{Discussion and Conclusions} \label{sec_Discu}
After having investigated each of the PDEs on its own, let us know discuss and draw conclusions from the results as a whole. 

Considering the solution time and accuracy, PINNs are not able to beat the finite element method in our study. Other than the inconclusive results for real and imaginary part in the Schr\"odinger 1D test, FEM solutions were generally faster at the same accuracy or at a higher accuracy.

After having solved the PDE, PINNs are sometimes faster at the pointwise evaluation of the respective solution: we were only able to show this in the 3D Poisson test. So when needing to evaluate a PDE on a very fine grid, one could consider solving a PINN. Although, in our examples, the solution time of FEM was so much faster that continued solving the PDE with FEM on different adapted grids would likely still be considerably faster than solving and evaluating the PINN.

We were particularly surprised that PINNs had difficulties with the Allen--Cahn equation using a  small $\varepsilon$. This might be due to the close-to-singular behaviour at the diffuse interface and the fact that PINNs aim to solve the possibly ill-conditioned strong form of the PDE. We anticipated that PINNs would outperform FEM:  the solution of an Allen--Cahn equation has very much the flavour of a classification (see \cite{budd}) at which neural networks excel; FEM requires a very fine grid to resolve the diffuse interface.  A similar case in which we were surprised that PINNs did not outperform FEM was the Schr\"odinger 2D examples, where the PDE solution, again, contains very finely structured areas. In both these cases an adaptive PINNs approach or variational PINNs \cite{VPINNs} might help. The latter would then also allow activation functions that have weak derivatives.

An aspect of the evaluation time that we have not considered throughout this work is the possibility of solving parameterised PDEs with neural networks, so-called operator approximators. See, for instance, Fourier Neural Operators \cite{Li2020} and DeepONets \cite{Lu2021}. Whilst the finite element method requires continued solutions of PDEs when changing the parameters, neural networks can take parameters as additional inputs and be trained throughout all of them. They have been shown to work well as surrogates  if the PDE needs to be solved sufficiently often; see also \cite{Maass2022}. In a future work, those should be compared to classical methods for parameterised PDEs, such as reduced bases \cite{Quarteroni2016}  or low-rank tensor methods \cite{kressner}. Again, the  time of the offline phase in which the parametric model is constructed or trained has to be considered carefully.

PINNs were good at the transition into higher dimensions: there is no increment in computational cost from the Poisson equation in 2D and 3D. This hints at the efficiency of PINNs in high-dimensional settings, in which classical techniques (such as FEM) are prohibitively expensive. This has been considered in \cite{kunisch}. In general, PINNs open up many interesting new research directions, especially when employed to solve such high-dimensional PDEs or when combining PDEs and data. The analysis of PINNs is both very challenging and highly interesting. Our study suggests that for certain classes of PDEs for which classical methods are applicable, PINNs are not able to outperform those. 

\begin{ack}
Initial research for this work was carried out with support from the Philippa Fawcett internship programme.

TGG and CBS acknowledge the support of the Cantab Capital Institute for the Mathematics of Information and the European Union Horizon 2020 research and innovation programme under the Marie Skodowska-Curie grant agreement No. 777826 NoMADS. TGG additionally acknowledges the support of the EPSRC National Productivity and Investment Fund grant Nr. EP/S515334/1 reference 2089694. JL and CBS acknowledge support from the EPSRC grant EP/S026045/1. CBS acknowledges support from the Philip Leverhulme Prize, the Royal Society Wolfson Fellowship, the EPSRC advanced career fellowship EP/V029428/1, EPSRC grants EP/T003553/1, EP/N014588/1, EP/T017961/1, the Wellcome Trust 215733/Z/19/Z and 221633/Z/20/Z, and the Alan Turing Institute.

We also acknowledge the support of NVIDIA Corporation with the donation of two Quadro P6000 GPU used for this research.
\end{ack}

\bibliographystyle{abbrv}
\bibliography{bibliography}
\end{document}